\documentclass[a4paper,oneside,10pt]{article}
\usepackage{amsmath,amsthm,amsfonts}
\usepackage[english]{babel}
\usepackage[latin1]{inputenc}
\usepackage[T1]{fontenc}
\usepackage{url}
\usepackage[text={16.5cm,24.5cm},centering]{geometry}
\usepackage{graphicx}

\newcommand{\corps}[1]{\mathbb{#1}}
\newcommand\R{\corps{R}}

\newcommand{\compl}[1]{#1^{c}}

\newcommand{\esp}{\mathbf{E}}
\newcommand{\var}{\mathbf{Var}}
\newcommand{\Ent}{\mathbf{Ent}}
\newcommand{\pseudoEnt}{\text{\textbf{Ps-Ent}}}

\DeclareMathOperator{\capacite}{Cap}
\newcommand{\grad}{\nabla}

\DeclareMathOperator{\osc}{osc}
\DeclareMathOperator{\argmin}{argmin}

\newcommand{\abs}[1]{|#1|}
\newcommand*{\nrj}[1]{\mathcal{E}(#1)}
\newcommand*{\Nrj}[2]{\mathcal{E}_{#1}(#2)}
\newcommand*{\rndNrm}{\mathcal{N}}
\newcommand*{\norme}[2]{\lVert#2\rVert_{#1}}
\newcommand*{\ninf}[1]{\norme{\infty}{#1}}
\let\orlicz=\norme
\newcommand*{\newOrlicz}[2]{%
  \expandafter\newcommand\csname orl#1\endcsname[1]{\orlicz{#2}{##1}}%
}

\newcommand{\ddt}[1]{\frac{d}{dt}\left(#1\right)}
\newcommand{\ddtUnSurSigma}{\ddt{\frac{1}{\sigma(t)}}}

\newcommand{\psiStar}{{\psi^\star}}
\newcommand{\phiStar}{{\phi^\star}}
\newOrlicz{Phi}{\phi}
\newOrlicz{Tau}{\tau}
\newOrlicz{Psi}{\psi}
\newOrlicz{PsiStar}{{\psi^\star}}
\newOrlicz{PhiStar}{{\phi^\star}}

\usepackage{amsthm}
\newtheorem{hpth}{Hypothesis}
\newtheorem{hpthU}{Hypothesis U}
\swapnumbers
\newtheorem{thm}{Theorem}
\newtheorem{prop}[thm]{Proposition}

\newtheorem{lmm}[thm]{Lemma}
\newtheorem*{rmq}{Remark}
\newtheorem{dfn}[thm]{Definition}
\newcommand{\ie}{\emph{i.e.}}
\newcommand{\eg}{\emph{e.g.}}

\newcommand{\cf}{\emph{cf.}}
\newcommand{\term}[1]{\emph{#1}} %
\newcommand{\ind}[1]{\mathbf{1}_{#1}}

\newcommand{\scal}[2]{\langle #1,#2\rangle}
\newcommand{\norm}{\abs}
\newcommand{\law}{\mathcal{L}}
\newcommand{\bra}[1]{\langle #1 \rangle}

\newcommand{\Dirac}[1]{\delta_{#1}}

\newcommand{\ball}[1]{\mathcal{B}_{#1}}
\DeclareMathOperator{\supp}{supp}
\newcommand{\hess}{\mathbf{Hess}}

\title{Annealing diffusions in a slowly growing potential}
\author{Pierre-André Zitt\thanks{Équipe Modal'X (EA3454), Université Paris \textsc{X}, Nanterre}}

\begin{document}
\maketitle

\begin{abstract}
  We consider a continuous analogue of the simulated annealing algorithm in $\R^d$, namely the solution of the SDE $dX_t = \sigma(t) dB_t - \grad V(X_t)dt$, where $V$ is a function called potential. We prove a convergence result, similar to the one in \cite{Mic91}, under weaker hypotheses on the potential function. In particular, we cover cases where the gradient of the potential goes to zero at infinity. The main idea is to replace the Poincaré and log-Sobolev inequalities used in \cite{Mic91,HCS87} by weak Poincaré inequalities (introduced in \cite{RW01}), and to estimate constants with measure-capacity criteria. We show that the convergence still holds for the ``classical'' schedule $\sigma(t) = c/\ln(t)$, where $c$ is bigger than a constant related to $V$.


\end{abstract}

\section*{Introduction}%

The goal of this article is to study a continuous analogue of a discrete optimization algorithm called \emph{simulated annealing}. This algorithm was introduced in 1983 by Kirkpatrick, Gelatt and Vecchi, and aims at  finding ``good'' (if not perfect) solutions to complex problems. The crucial idea is to perturb the standard gradient descent by a random noise; hopefully this noise will get the process out of traps (local minimas), and help it reach the global minimum. The noise is taken relatively big at the beginning, so that the process explores the space, and is gradually reduced thereafter.

The standard case is the discrete case (in time and space); here we consider a process on $\R^d$ in continuous time. Note that more complicated state spaces have been studied, see for example \cite{Jac94,JR95,Jac96}; here we will stick to $\R^d$. This ``annealing diffusion'' process has already been studied by several authors. Hwang, Chiang and Sheu (\cite{HCS87}) proved the convergence under quite strong assumptions, using comparisons with the associated (ordinary) differential equation and results on the trajectories (estimates of exit times from domains, etc.). The result was enhanced by Royer (\cite{Roy89}). The approach we follow was developed by L.~Miclo in \cite{Mic92} (and in his doctoral dissertation \cite{Mic91}), and reduces the problem to the convergence of a single quantity, the free energy. Since then, other questions have been asked: speed of convergence, choice of a better algorithm etc. (see \eg{} the survey \cite{Loc00}). Let us also note that the ``functional inequalities'' approach has also been used extensively for other (possibly discrete) models, and other closely related algorithms (see \eg{} \cite{DMM99} for a study of a generalized simulated annealing process).

A common feature of these works on global optimization on $\R^d$ is the quite strong assumptions they require on the growth of the potential. In particular, the norm of the gradient is supposed to go to infinity at infinity. These hypotheses are technically useful: they guarantee that, at any fixed temperature, the generator has a spectral gap, which in turn gives estimates on the rate of convergence. Let us note that the ``cooling schedule'' (\ie{} the choice of the temperature as a function of time) for which the process converges is linked with the speed of explosion of the  spectral gap, but that it can be read directly on the potential (see below the remarks on the constant $d^\star$).

A natural question arises: what happens when the gradient of the potential does not go to infinity, and when there is no spectral gap? Do we need to change the cooling schedule to reflect the slow-down of the diffusions at fixed temperature, or does the local structure of the potential dictate the optimal schedule?

Before we answer this question, let us be more precise and give our hypotheses.

\bigskip

We study the following optimization problem: how to find the minimum of a function $V$ on the space $\R^d$. To solve this problem, we introduce the following stochastic differential equation:
\[
\left\{
\begin{array}{rl}
    dX_t& = \sqrt{\sigma(t)}dB_t - \frac{1}{2}\nabla V(X_t)dt, \\
    X_0 &\sim m_0.
  \end{array}
\right.
\]
The function $\sigma$ will be called \emph{temperature}, and will be a (deterministic) function of time, decreasing to zero.

Intuitively, this process is similar to simulated annealing: we perturb a gradient descent by a stochastic term whose intensity decreases over time. 

We would like to know if the process finds a point where the global minimum is reached; we will show that it does, in a weak sense. 
\begin{dfn}
  The annealing process starting from a law $m_0$ is said to converge if its law $m_t$ at time $t$ converges weakly to a measure supported by $\argmin V$. In particular, if the global minimum of $V$ is reached in a single point $x_0$, the process converges if $m_t$ goes to a Dirac mass at $x_0$.
\end{dfn}

Let us now recall the result we would like to generalize: this is the main result of \cite{HCS87,Roy89,Mic92},  as it appears in \cite{Mic92}.
\begin{thm}[L.~Miclo]
  \label{thm=Miclo}
  If $V$ satisfies some regularity assumptions, and the following conditions:
  \begin{itemize}
    \item 
   $V \xrightarrow[x\to\infty]{} \infty, $
 \item 
   $\abs{\grad V} \xrightarrow[x\to\infty]{} \infty, $
 \item
   $\abs{\grad V} - \Delta V \text{ is bounded from below},$
   \end{itemize}
   then there exists a constant $d^\star$ such that, for any $c>d^\star$, and for $\sigma(t) = c/\ln(t)$, the annealing process converges. 
\end{thm}
To understand the direction in which we generalize this result, let us note that this theorem applies for any potential $V$ which is equal to $\abs{x}^\alpha$ outside a compact set, whenever $\alpha$ is strictly bigger than $1$. It is then a quite natural question to ask whether this still holds when $\alpha$ is strictly less than $1$. Our hypotheses, which we now state, allow us to treat this case.

\newcommand{\croissanceV}{{m_V}}

\begin{hpth}[Global minimum]%
  \label{hpth=globalMinimum}
  The potential has a unique global minimum, located at the origin and $V(0) = 0$. Moreover, this minimum is non degenerate: $\hess V (0)$ is positive definite.
\end{hpth}%

\begin{hpth}[Growth at infinity]%
  \label{hpth=croissanceDeV}
  The potential $V$ goes to infinity at infinity faster than a logarithm: 
  \[ \exists \croissanceV >1,\exists C, \quad 
  V(x) \geq \ln(\abs{x})^{\croissanceV } - C.\]
\end{hpth}%

\begin{hpth}[Bounded gradient]%
  \label{hpth=gradientBorne}
  The potential $V$ is continuously differentiable, and its gradient is bounded:
  \[ \ninf{\grad V} < \infty.\]
\end{hpth}%

\begin{hpth}[Concavity]
  \label{hpth=concavite}
  The Laplacian of $V$ is negative at infinity: there exists a compact set $K$ compact such that
  \[ \forall x\notin K, \quad \Delta V(x) \leq 0.\]
\end{hpth}

One last hypothesis will be added in section \ref{sec=multidim}, regarding the structure of local minima of $V$.

\bigskip

These hypotheses call for a few remarks.

The first one simplifies the problem at hand: there is only one goal to go after. If the weak limit of the equilibrium measures $\mu_\sigma$ (\cf{} infra) is known (some results in this direction may be found in \cite{Mic92,Hwa82}), the arguments given here should work in the same way. The non-degeneracy hypothesis may be weakened too (see \eg{} section \ref{sec=unidim} for a slight generalization in $d=1$)
However, this restriction allows for two simplifications: it gives an estimate of the partition function $Z_\sigma$, and avoids more intricate reasonings in the computation of the weak inequalities (section \ref{sec=multidim}).

The growth hypothesis is not very restrictive. In particular, $V$ may grow like $\abs{x}^\alpha$ with $\alpha<1$ (or even slower). These cases were not covered in the literature. Let us note that we do not know what happens in the limit case (when $\croissanceV = 1$, \ie{} the tails of the equilibrium measures are polynomial).

In the light  of the previously known results, the bounded gradient assumption seems less stringent: in some sense, we already know what happens when the gradient is big. The hypothesis could probably be lifted if we allowed a polynomial growth, or a control by $V$, but we keep it for the sake of clarity. 

Finally, the condition on $\Delta V$ seems more restrictive. It will only be used in the proof of the moment bound (section \ref{sec=controleDeMoments}). It could probably be replaced by a condition like $\Delta V \leq C \abs{\grad V}^2$. However, in the ``natural example'' where $V(x) = \abs{x}^\alpha$ at infinity, the Laplacian is indeed negative  if $\alpha<1$, and this example was one motivation for investigating the problem. 
Moreover, even this weakened hypothesis would not allow the existence of traps at infinity, however shallow they may be. It would be interesting to know what could happen if there were such traps: either they have no effect (in the sense that the same cooling schedule may be chosen), or they slow down the process too much and destroy the convergence.

\bigskip

Our principal result is the following.
\begin{thm}
  If the potential $V$ satisfies the hypotheses above, there exists a constant $d^\star$ such that,
  if we choose
  \[ \sigma(t)  = \frac{c}{\ln(t)}, \]
  with $c>d^\star$, the annealing process converges.
\end{thm}
This result generalizes theorem \ref{thm=Miclo} by allowing more general choices for the potential function. In particular, as we will see in the sequel, the equilibrium measures need not satisfy a Poincaré inequality. Nonetheless, the critical cooling schedule is the same, which contradicts the intuition that the speed was given by   the Poincaré constants. In fact, what seems to prevail is the behavior of $V$ in a compact set, and from a certain point of of view, that is precisely what the weak inequalities capture.

The remainder of the paper is organized in the following way. Firstly, we explain the analytic approach of L.~Miclo and give the main line of the proof. 

This proof, under our weakened hypotheses, uses \emph{weak Poincaré inequalities}. We will need controls over their dependence on temperature: these are established in sections \ref{sec=unidim} and \ref{sec=multidim}, respectively in the one- and multi-dimensional case. These three sections are the core of the proof of the convergence result. 

The quite technical \ref{sec=CapMesEtPoincareFaible}\textsuperscript{th} section gathers definitions and results about Orlicz norms and weak inequalities. Finally, we postpone to the annexes a comparison between functions centered by their mean or by their median, a  moment  bound for the annealing process, and a brief proof of the estimation of the partition function.

\section{The convergence of the process (main line of the proof)} %
\subsection{A differential inequality for  the free energy} %
Before we describe the main idea, we introduce some notation.
Consider the SDE defining the annealing diffusion, but with a constant temperature $\sigma$. The process is then  a classical diffusion with a gradient drift. The corresponding generator is given by:
\[
  L_\sigma : f \mapsto \frac{\sigma}{2} \Delta f- \frac{1}{2}\nabla V \nabla f.
\]
The measure $\mu_\sigma$ defined by
\[d\mu_\sigma = \frac{1}{Z_\sigma} \exp\left( - \frac{V}{\sigma} \right) d\lambda,\]
is reversible for this process ($Z_\sigma$ is a normalization constant).
We will call $\mu_\sigma$ the \term{instantaneous equilibrium measure}.

It's easy to see that, as $\sigma$ goes to zero, the measures $\mu_\sigma$ concentrate around the global minimum of the potential (which is found at the origin by hypothesis). In fact, we even have the following convergence.
\begin{prop}
   The measures $\mu_\sigma$ converge weakly: 
  \[ \mu_\sigma \xrightarrow[\sigma\to 0]{} \Dirac{0}.\]
  Moreover, the normalization constant $Z_\sigma$ behaves like $\sigma^{d/2}$.
\end{prop}
The asymptotic behavior of $Z_\sigma$ is proved in annex \ref{sec=fonctionDePartition}. 

\bigskip
In order to prove the convergence of the process, we follow the approach of L.~Miclo (\cite{Mic92}) and show that the relative entropy of the law of the process with respect to its instantaneous equilibrium measure goes to zero.

More precisely, let $f_t$ be the density of $m_t = \mathcal{L}(X_t)$ with respect to the equilibrium measure $\mu_t$. The relative entropy
(also called free energy) is $I_t = \int f_t \log f_t d\mu_t$, which can be rewritten as $I_t = \Ent_{\mu_t}(\sqrt{f_t}^2)$. The finiteness of $I_t$ is established in annex \ref{sec=finiteEntropy}.
We would like to study the evolution of $I_t$; the natural idea is to differentiate it. One can justify the following formal computation:
 \begin{prop}[Differentiation of the free energy]%
   The derivative of the free energy is given by:
  \begin{align*}
  \frac{dI_t}{dt} &= \frac{1}{\sigma(t)^2}\sigma'(t) \int
  V\times(1-f_t) d\mu_t - 2 \sigma(t)\int\abs{\nabla\sqrt{f_t}}^2 \\
                  &= \frac{1}{\sigma(t)^2}\sigma'(t) \int
  V\times(1-f_t) d\mu_t - 2 \sigma(t)\Nrj{\mu_t}{\sqrt{f_t}}.
\end{align*}
\end{prop}%
\begin{rmq}%
By $\Nrj{\mu t}{f}$ we denote $\int \abs{\nabla f}^2 d\mu_t$. This is somewhat improper --- strictly speaking, this is the energy associated
with the generator $(1/2)\Delta - (1/(2\sigma))\nabla V \nabla$ (so we should multiply our energy by $\sigma$ to get the ``real'' one). However, 
the classical criteria for functional inequalities are written for this form of the energy. 
\end{rmq}%

The first term is set aside for the time being, we shall bound it later directly by a function of $t$.

Following the classical path leading from functional inequalities to semigroup estimates, we now try to control the energy term on the right hand side.

If the measures $\mu_t$ satisfied logarithmic Sobolev inequalities, everything would be fine: the energy of $\sqrt{f_t}$ could be controlled by its entropy with respect to $\mu_t$, and we would get $I_t$ back on the right hand side of the inequality. We would still have to know how the constants in the logarithmic Sobolev inequality depend on the small parameter $\sigma$, and get an upper bound for the first term, but we could get the convergence of $I_t$ to zero.

Unfortunately, the scaling behavior of the constants in the logarithmic Sobolev inequality (\ie{} the way they behave when $\sigma$ goes to zero) is not clear.
Moreover,  this inequality need not  hold, and in fact it won't under our hypotheses.

In Miclo's paper, the first difficulty is overcome thanks to  a Poincaré inequality, weaker than the logarithmic Sobolev inequality, but for which the constants are well known. However, even this  inequality won't be satisfied in our case, and we have to find another way. 

Our idea is to consider a still weaker functional inequality, namely a \term{weak Poincaré inequality}, written with an Orlicz norm. Weak Poincaré inequalities were introduced by M.~Röckner and F.-Y.~Wang in \cite{RW01}, originally with an $L^\infty$ norm and the mean of $f$ instead of a median on the right hand side. We will give a brief account on weak inequalities and Orlicz norms in section \ref{sec=CapMesEtPoincareFaible}, and explain the link between the original inequality and the one we use. 

For now, let us just state this inequality. It reads: 
\begin{equation}
  \label{eq=defPoincareFaible}
  \forall f,\forall r,\qquad 
  \var_{\mu_t}(f) \leq \alpha_t(r)\Nrj{\mu_t}{f}
    + r \orlPhi{f - m_f}^2,
\end{equation}
where $m_f$ is a median of $f$ under $\mu$, and $\alpha_t$, a decreasing function of $r$, is the compensating function.
The Orlicz norm is not easily tractable, but we will see (\cf{} lemma \ref{lmm=controleOrliczEntropie}) that it can be bounded by the entropy: there exists a $C$ such that, for all positive $f$,
\[\orlPhi{f - m_f}^2 \leq C(\mu(f^2)+\Ent(f^2)).\]
At this point, the energy is bounded above by three terms: $\mu(f^2)$, the entropy of $f$ and its variance. 
To get rid of the variance term, we would like to bound it by entropy-like quantities. To this end we introduce the following definition.
\begin{dfn}%
  For any probability measure $\mu$ and any positive $f$, we will call \term{pseudo-entropy} the quantity:
  \[\pseudoEnt(f) = \int f \log^2\left(e+
  \frac{f}{\norme{1,\mu}{f}}\right)d\mu.\]
\end{dfn}%
With this definition in hand, we can state (\cite{Mic92}, lemma 4):
\begin{lmm}%
  There exists a $\delta_0$ such that, for all probability measure $\mu$ and all positive $f$ with $\mu(f^2)=1$, 
  \[\forall \delta<\delta_0, \frac{1}{\delta}\var_\mu(f) + 4
  \delta\pseudoEnt_\mu(f^2) \geq \Ent_\mu(f^2).\]
\end{lmm}%
Let us put all these inequalities together: we get that for all probability measure $\mu$, if $\mu$ satisfies the 
weak Poincaré inequality \ref{eq=defPoincareFaible}, then for all positive $f$ with $\int f^2 d\mu = 1$, 
\[
  \delta \Ent_\mu(f^2) - 4 \delta^2 \pseudoEnt(f^2)%
  \leq  \var_\mu(f) \leq \alpha(r)\Nrj{\mu}{f,f}
    + Cr\Ent_\mu(f^2) + Cr.
\]
This entails a lower bound on the energy: 
\[
\Nrj{\mu}{f,f} \geq - \frac{4}{\alpha(r)}\delta^2 \pseudoEnt(f^2) %
  - C\frac{r}{\alpha(r)}  %
  + \frac{1}{\alpha(r)}(\delta - Cr)\Ent_\mu(f^2).
\]
Let us get back into our special case, and take $\mu=\mu_t$, $f = \sqrt{f_t}$. The entropy $\Ent(f^2)$ just becomes $I_t$,
and we can plug the inequality back in the differential equation for $I_t$:
\newcommand{\circledNumber}[1]{\text{\textcircled{\scriptsize #1}}}

\newcommand{\Un}{\circledNumber{1}}%
\newcommand{\Deux}{\circledNumber{2}}%
\newcommand{\Trois}{\circledNumber{3}}%
\newcommand{\Quatre}{\circledNumber{4}}%
\begin{align*}%
\frac{dI_t}{dt}%
 & \leq  \frac{1}{\sigma(t)^2}\sigma'(t) \int  V\times(1-f_t) d\mu_t%
  + 8\delta^2\frac{\sigma(t)}{\alpha_t(r)} \pseudoEnt(f_t)
  + 2C\sigma(t) \frac{r}{\alpha_t(r)}%
  - 2(\delta-Cr)\frac{\sigma(t)}{\alpha_t(r)} I_t
\end{align*}%
Since $\sigma$ is non-increasing in time, we may omit the $1$ in $(1-f_t)$ in the first term, and since $f_td\mu_t = dm_t$, 
\begin{align}
  \label{eq=inequationDiffPourIt}
\frac{dI_t}{dt}%
 & \leq  \ddtUnSurSigma \int  V dm_t%
  + 8\delta^2\frac{\sigma(t)}{\alpha_t(r)} \pseudoEnt(f_t)
  + 2C\sigma(t) \frac{r}{\alpha_t(r)}%
  - 2(\delta-Cr)\frac{\sigma(t)}{\alpha_t(r)} I_t
\end{align}
Our goal is to obtain a differential inequality involving only $I_t$ and explicit functions of $t$, so that we may deduce information on the evolution of $I_t$. Since $\sigma$ is known, this leaves us with three questions. First, we have to obtain controls on $\int V dm_t$ and on the pseudo-entropy --- we will get explicit bounds in $t$. Once this is done, we have to estimate the compensating function $\alpha_t$. Finally we must choose $r$ and $\delta$ depending on $t$ in a suitable way, so that the inequality on $I_t$ is good enough to prove the convergence to zero.

We now deal with the first problem.

\subsection{Moment bounds and pseudo-entropy} %
The first inequality is a moment bound on the value of the potential at time $t$. The proof is postponed to the annexes.
\begin{lmm}\label{lmm=controleV} 
 Suppose that hypotheses \ref{hpth=gradientBorne} and \ref{hpth=concavite} hold, and that the initial law  $m_0$  satisfies: $\int V^p m_0(dx) <\infty$. Then there exists an $M$ such that:
 \[\int V^p(x)m_t(dx) \leq M\sigma(t)^p \ln(t)^p (\ln\ln t)^{3p}.\]
\end{lmm}
The last result will be used directly, but it also helps us prove the following bound.
\begin{lmm}\label{lmm=controlePseudoEntropie}
Suppose that $\int V^2dm_0$ is finite, and that the cooling schedule has the form: $\sigma(t) = c/\ln(t)$, for a positive constant $c$. Then there exists an $A$ such that, for all big enough $t$, 
\[\pseudoEnt(f_t) \leq A\ln(t)^{2} (\ln\ln(t))^6.\]
\end{lmm}
\begin{proof}%
We differentiate the quantity under scrutiny, namely $J_t = \pseudoEnt_{\mu_t}(f_t)$. The following formal computation can be justified (\cf{} \cite{Mic92}):
\[
\frac{dJ_t}{dt} =%
  - \frac{\sigma(t)}{2} \int F'(f_t) \abs{\nabla f_t}^2 d\mu_t%
  + 2 \ddtUnSurSigma \int\log (e+f_t) \frac{f_t}{e+f_t} \left(%
     V - \int V(x)dm_t
   \right)dm_t,
\]
where $F(x)  =\frac{2x}{x+e}\log(x+e) + \log^2(x+e).$ Since $F$ is non decreasing (in $x$), and $\sigma$ is 
positive, the first term is bounded above by $0$. Moreover, since $V$ is positive and $1/\sigma$ increases, 
we may also forget the $\int V(x) dm_t$ in the second term. We get: 
\begin{align*}%
  \frac{dJ_t}{dt}%
  &\leq 2\ddtUnSurSigma\int
    \frac{f_t}{e+f_t}\log(e+f_t) V dm_t\\
  &\leq 2\ddtUnSurSigma\int \log(e+f_t)Vdm_t\\
  &\leq 2\ddtUnSurSigma%
    \left(\int \log^2(e+f_t) dm_t\right)^{\frac{1}{2}}%
    \left(\int V^2dm_t\right)^{\frac{1}{2}}\\
&= 2\ddtUnSurSigma J_t^{\frac{1}{2}}%
 \left(\int V^2dm_t\right)^{\frac{1}{2}}.
\end{align*}%
After dividing by $2J_t^{\frac{1}{2}}$, the left hand side becomes the derivative of $\sqrt{J_t}$. The right hand
side may then be bounded (\cf{} previous lemma):
\begin{align*} %
\frac{d\sqrt{J_t}}{dt}
&\leq \ddtUnSurSigma\left(\int V^2dm_t\right)^{\frac{1}{2}}\\
&\leq \ddtUnSurSigma \sqrt{M}
  \sigma(t) \ln(t)) (\ln\ln(t))^3.
\end{align*}
The explicit value of $\sigma$ allows us to simplify: 
\begin{align*}
\frac{d\sqrt{J_t}}{dt}
&\leq \sqrt{M} \frac{1}{t} \left(\ln\ln(t)\right)^3.
\end{align*}%
An easy computation shows that the right hand side may be bounded by:
\[\sqrt{M} \frac{d}{dt}\left( \ln(t)(\ln\ln(t))^3\right).\]
To conclude the proof, we integrate this inequality between a (fixed and big enough) $t_0$ and the current time $t$. 
The constant $A$ naturally depends on the initial law $m_0$ (through the value of $M$ and through the pseudo-entropy at time $t_0$). 
\end{proof}%

\subsection{From the differential inequality to the convergence of the entropy} %
It is now time to get back to our differential inequality and apply the bounds we just derived.
We fix a logarithmic cooling schedule:
\[ \sigma(t) = \frac{c}{\ln(t)}.\]
Recall that we showed (inequality \ref{eq=inequationDiffPourIt}):
\begin{align*}
\frac{dI_t}{dt}%
 & \leq  \ddtUnSurSigma \int  V dm_t%
  + 8\delta^2\frac{\sigma(t)}{\alpha_t(r)} \pseudoEnt(f_t)
  + 2C\sigma(t) \frac{r}{\alpha_t(r)}%
  - 2(\delta-Cr)\frac{\sigma(t)}{\alpha_t(r)} I_t
\end{align*}
We use the moment bound (lemma \ref{lmm=controleV}) to deal with the first term, and lemma \ref{lmm=controlePseudoEntropie}
to bound the second one.
\begin{align*}
\frac{dI_t}{dt}%
 & \leq  \ddtUnSurSigma M \ln\ln(t)^3
  + 8M\delta^2\frac{\sigma(t)}{\alpha_t(r)}
  (\ln(t))^{2} (\ln\ln(t))^6
  + 2C\sigma(t) \frac{r}{\alpha_t(r)}%
  - 2(\delta-Cr)\frac{\sigma(t)}{\alpha_t(r)} I_t  \notag\\
\end{align*}
We number our four terms and define:
\begin{align*}
 \Un &=  \ddtUnSurSigma M \ln\ln(t)^3
 &
 \Trois &=  2C\sigma(t) \frac{r}{\alpha_t(r)} \\
 \Deux &=  8M\delta^2\frac{\sigma(t)}{\alpha_t(r)}
  (\ln(t))^{2} (\ln\ln(t))^6 
 &
 \Quatre &=  2(\delta-Cr)\frac{\sigma(t)}{\alpha_t(r)}\\
\end{align*}
The inequality becomes:
\begin{equation}
  \frac{dI_t}{dt} = \Un + \Deux + \Trois  - \Quatre I_t \label{eq=UnDeuxTroisQuatre}
\end{equation}
This last inequality will allow us to prove that the free energy goes to zero. To this end, we use the same lemma as L.~Miclo:
\begin{lmm}\label{lmm=equaDiffClassique}%
Let $I$ be a positive function, and suppose:
\[\frac{dI_t}{dt}\leq a(t) - b(t)I(t),\]
where $a,b$ are positive functions and satisfy:
\begin{enumerate}
\item $\int^\infty b(t) = \infty,$
\item $\frac{a(t)}{b(t)}\xrightarrow{t\to\infty} 0.$
\end{enumerate}
Then $I$ goes to zero when $t$ goes to infinity.
\end{lmm}%
Our goal is now to use the inequality \ref{eq=UnDeuxTroisQuatre} to check the hypotheses of this lemma. 
We choose $\delta$ and $r$ as follows.
\begin{equation}
  \label{eq=choixDeR}
  \begin{cases}
  \delta_t &= \frac{1}{\ln(t)^2 (\ln\ln(t))^7}\\
  r_t   &= \frac{1}{ C \ln(t)^2 (\ln\ln(t))^8}
\end{cases}
\end{equation}
This ensures:
\begin{align*}
\frac{\Trois}{\Quatre}  &= \frac{Cr_t}{\delta_t - Cr_t}
  \sim \frac{C}{\ln\ln(t)} \rightarrow 0,\\
\frac{\Deux}{\Quatre} &= \frac{4M\delta_t^2}{\delta_t - Cr_t} \ln^2 (t) (\ln\ln(t))^6
  \sim 4M \delta_t \ln^2( t) (\ln\ln(t))^6
\rightarrow 0.
\end{align*}
Two things remain to check:
\begin{align*}
  \frac{\Un}{\Quatre}&\rightarrow 0  & &\text{and} &
  \int^\infty\Quatre &= \infty.
\end{align*}
This is where we need bounds on the weak Poincaré inequalities: we have to know how $\alpha_t$ behaves for our particular choice of $r$. This is the aim of the following sections, in one or many dimensions. 

In both cases, we will get:
\begin{lmm}
  \label{lmm=majAlpha}
  There exists a constant $d^\star$ such that, for all    $D^\star>d^\star$, 
  \[
  \exists C_\alpha \qquad \alpha_t(r_t) \leq C_\alpha \exp\left( \frac{D^\star}{\sigma} \right).
  \]
  For the cooling schedule $\sigma(t) = c/\ln(t)$, we get:
  \[
   \alpha_t(r_t) \leq C_\alpha t^{D^\star/c}.
  \]
\end{lmm}
In the one-dimensional case, this follows from theorem \ref{thm=alphaUnidim} below, and the choice of $r_t$. The multi-dimensional case is proved in theorem \ref{thm=PoincareDefectif} and the discussion that follows it. 
\begin{rmq}
  The approach in the one- and multi-dimensional case will differ slightly. In the former, we prove a (full) weak Poincaré inequality, \ie{} we estimate the whole function $\alpha_t$, and then use this estimate at the point $r_t$. In the latter, we will only prove a bound on $\alpha_t$ at $r_t$ and disregard the other points.
\end{rmq}

\bigskip

We may know get back to our proof. Recall that we have assumed:
\[ \sigma(t) = \frac{c}{\ln(t)}, \qquad c > d^\star,\]
so that we may always pick a $D^\star$ strictly less than $C$.

Let us check the two remaining points. First we must prove that $\Un/\Quatre$ converges. Since $\sigma$ is explicit and we know a bound on  $\alpha(r)$, we see that:
\begin{align*} %
  \frac{\Un}{\Quatre} 
  & = \ddtUnSurSigma M (\ln\ln(t))^3 \times \frac{%
    \alpha_t(r_t) 
    }{
    2(\delta - Cr)\sigma(t)
    } \\
    &\leq M'\frac{1}{t} \ln(t)^3\ln\ln(t))^{10}  \alpha_t(r_t).
 \end{align*} %
 where $M,M'$ are constants.

Using the bound on $\alpha$ we just recalled (lemma \ref{lmm=majAlpha}), we get:
 \begin{align*}
   \frac{\Un}{\Quatre}
   &\leq 
   M'' \frac{t^{D^\star/c}}{t} \left( (\ln t)^3(\ln\ln t)^{10}\right)\\
\end{align*}
Since $c>D^\star$, $\Un/\Quatre$ goes to zero, as was claimed.

Just in the same way, we have, for $t$ big enough:
\begin{align*}
  \Quatre  &= 2(\delta_t - Cr_t)\frac{\sigma(t)}{\alpha_t(r_t)} \\
  &\sim M'' (\ln t)^{-3} (\ln\ln t)^{-7} \frac{1}{t^{\frac{D^\star}{c}}}.
\end{align*}
Once more, the condition $c>D^\star$ guarantees that the integral of this quantity diverges, which was expected.

This allows us to apply lemma \ref{lmm=equaDiffClassique}, and prove that $I_t$ converges to $0$. Thanks to Pinsker's inequality, the total variation between $m_t$ (law of the process) and $\mu_t$ (the instantaneous equilibrium) converges too. Since we already know that $\mu_t$ converges weakly to the Dirac mass $\delta_0$, this concludes the proof.

\subsection{Some remarks} %
Our theorem immediately raises a few questions. Some of these  have already been asked when we discussed the hypotheses --- equilibrium measures with polynomial tails are not covered, and we do not know what happens when there are traps at infinity.

It would also be interesting to know what happens if we cool faster than the ``good'' schedule. \emph{A priori}, the process has no reason to converge to the global minimum; intuitively it should freeze in some local trap. One could ask if this trap is a good approximation of the global aim. Answering this question seems impossible in all generality, one should have to assume much more on the potential function, and on the starting point. The ``analytic'' approach may not be the best suited for this task.

\section{The one-dimensional case}%
\label{sec=unidim}
\newcommand{\dStar}{{d^\star}}
\newcommand{\deltaStar}{{\delta^\star}}
In this section we treat the case of a one-dimensional potential, for which we derive a weak Poincaré inequality (more precisely we prove lemma \ref{lmm=majAlpha}). 

The major advantage of this case is that, in one dimension, explicit (Hardy-like) criteria are known for weak inequalities. Thus we are able to prove a quite general result (the de-coupling of the parameters $s$ and $\sigma$ in the weak inequality). This has a small price: we restrict ourselves to potentials that grow like a power of $x$, and do not cover the case $V(x) = \log(x)^\alpha$ at infinity (for some $\alpha>1$). It should be noted that the multidimensional argument (\cf{} next section) may still be used in this logarithmic case.

Let us write down a few notations. The potential $V$ is a real function, continuously differentiable. For any (small) $\sigma$, we denote by $V_\sigma$ the function $\frac{1}{\sigma}V$, and by  $Z_\sigma = \int e^{-V_\sigma(x)}dx $ the partition function. We normalize $V_\sigma$ by defining $\Phi_\sigma$: $ \Phi_\sigma = V_\sigma + \log Z_\sigma$. 
The equilibrium measure $\mu_\sigma$ reads:
    \[ d\mu_\sigma = \frac{1}{Z_\sigma}\exp( - V_\sigma) d\lambda =  \exp( - \Phi_\sigma ) d\lambda. \]

    We now state our hypotheses on $V$. We suppose there exists a compact set $[K_1,K_2]$ such that the following holds.
\begin{hpthU}[Behavior near the minimum]
  \label{hpth=VPresDuMin}
  In $[K_1,K_2]$, $V$ is bounded below by $0$ and  above $V(K_1) = V(K_2)$. It reaches its minimum only once, at $x_1$. Near this point, $V$ behaves like: 
    \[V(x) \sim (x-x_1)^b,\]
  with $b>1$. Finally, there exists $\delta$ such that $V$ is bijective from $[x_1,x_1 + \delta]$ onto its image, and from $[x_1 -\delta,x_1]$ onto its image.
\end{hpthU}
This generalizes a little the general assumptions on the minimum: if $\hess V$ is positive definite at $x_1$, it satisfies 
this hypothesis with $b=2$.
  
\begin{hpthU}[Behavior outside the compact]
  \label{hpth=VLoinDuMin}
  Outside the compact, $V'$ and $\abs{V''}/(V'^2)$ are bounded:
  \begin{equation}%
    \label{equ=majVseconde}
    \exists C_V \forall x\notin [K_1,K_2],\qquad \frac{\abs{V''}}{V'^2} \leq C_V.
  \end{equation}
  In particular, $V'$ has no zero, $V$  decreases before $K_1$ and increases after $K_2$.
\end{hpthU}
\begin{hpthU}[The function $\beta$]
  \label{hpth=beta}
  There exists a function $\beta$ such that, for all $x$ outside the compact, 
  \begin{equation}
  \label{equ=beta}
    \beta \left(%
      \frac{\exp (-V(x))}{V'(x)}
    \right)
    \geq \frac{1}{V'(x)^2}.
  \end{equation}
\end{hpthU}
To apply the result to the annealing diffusion, we need an additional growth condition on $\beta$:
\begin{hpthU}[Behavior of $\beta$ near the origin]
  \label{hpth=croissanceBeta}
  There exist constants $A,C$ such that, near $0$, the following holds:
  \[\beta(s) \leq C \left(\log(\frac{1}{s})\right)^A.\]
\end{hpthU}
\begin{rmq}
We shall note here that the last two hypotheses hold if $V(x) = x^\alpha$ outside  a compact, with $\alpha \in (0,1]$, if we choose $\beta = C \left( \log(1/s) \right)^{\frac{4}{\alpha} -4}$ (\cf{} \cite{RW01,BCR05}). 
If $v$ grows like a logarithm to some power, this is not true ($\beta$ behaves like a power of $s$).  This explains the small loss of generality we spoke about above.
\end{rmq}

 We define, for all $x\geq x_1$,  $i(x) = \inf\{V(y),y>x\}$ and $s(x) = \sup\{ V(y),y\in[x_1,x]\}$. In the same way, $i(x) = \inf\{V(y),y<x\}$ and $s(x) = \sup\{V(y),y\in[x,x_1]\}$ for $x$ less than  $x_1$.  
    
Outside $[K_1,K_2]$, we have $i=V=s$, so $s-i$ is continuous with compact support. We call $\dStar$ its maximum value.

\begin{figure} %
  \includegraphics[width=.9\linewidth]{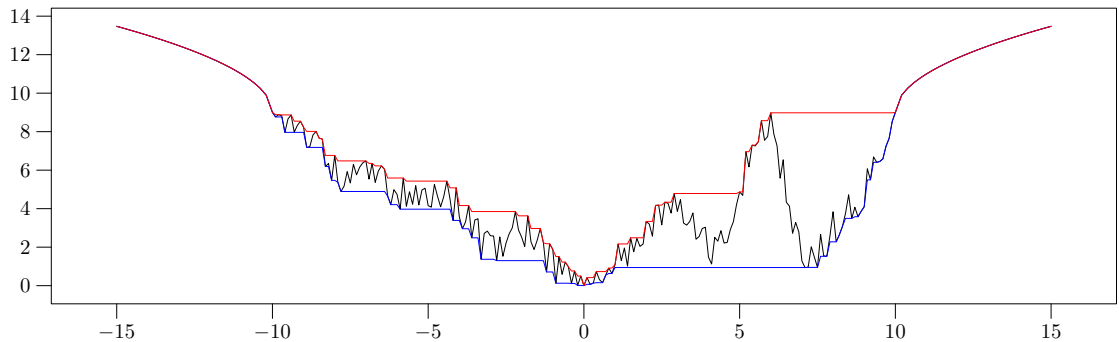}
  \caption{The potential $V$ and the associated functions $i$ and $s$. Here, $x_1 = 0$, $K_1 = -10$, $K_2 = 10$.  }
  \label{fig=doublePuits}
\end{figure} %

\bigskip

The main result of this section may now be stated as follows.
\begin{thm} \label{thm=alphaUnidim} %
  The measure $\mu_\sigma$ satisfies a weak Poincaré inequality with the $L^\infty$ norm, with a compensation 
  function $\beta_\sigma$ defined by:
  \begin{align*}
    \beta_\sigma(s )  &= C \exp\left( \frac{\dStar}{\sigma}\right) \beta(s),
  \end{align*}
  where $\beta$ is given by the hypothesis. Similarly, $\mu_\sigma$ satisfies a weak inequality with an Orlicz norm and the modified function $\alpha_\sigma$ given by:
  \[
  \alpha_\sigma(r) = C\beta_\sigma\left(%
    C'\exp(-\frac{4}{r})
    \right)
    = C \exp\left(\frac{\dStar}{\sigma}\right)
      \beta\left(
        C' \exp\left(- \frac{4}{r} \right)
      \right).
  \]
  Finally, there exists a constant $A$ such that the following bound holds:
  \[\alpha_\sigma(r) 
  \leq C \exp\left(\frac{\dStar}{\sigma}\right) \frac{1}{r^A}.
  \]
\end{thm} %

To prove this, we will use a result from Barthe, Cattiaux and Roberto (\cite{BCR05}, theorem 3), which gives estimates on the compensating functions for the $L^\infty$ norm. We will then use capacity-measure criteria to derive the result with the Orlicz norm. To state the result we need, we first give some additional notation. 

    Let $m_\sigma$ will be a median of $\mu_\sigma$, and for all $x$, 
\begin{align} %
  \notag 
  B_\sigma(x)%
  &= \frac{
      \int_{m_\sigma}^x e^{\Phi_\sigma(y)} dy%
      \times
      \int_x^\infty e^{ - \Phi_\sigma} dy
     }{
       \beta\left(\int_x^\infty e^{-\Phi_\sigma(y) }dy \right)
     }\\
\label{eq=defBSigma}
 B_\sigma%
 &= \sup_{x\geq m_\sigma} B_\sigma(x).
\end{align} %
By symmetry, we also define $b_\sigma(x)$ and $b_\sigma$ for $x\leq m_\sigma$.

The result from \cite{BCR05} reads:
\begin{thm}
  \label{thm=CapMesBCR}
  Let $\beta: (0,1)\to \R^+$ be non increasing, and $B_\sigma,b_\sigma$ be defined by \eqref{eq=defBSigma}.

  Then $\mu_\sigma$ satisfies the following weak Poincaré inequality~:
  \[ \var_{\mu_\sigma}(f) \leq C_\sigma \beta(s) \int \abs{\grad f}^2d\mu_\sigma + d \osc(f)^2,\]
  where $C \leq 12\max(b_\sigma,B_\sigma)$.
\end{thm}
Note that their result is actually  stronger, since it also gives a lower bound on the optimal constant $C$ in terms of some quantities very similar to $B_\sigma$.

\bigskip

To use this result, we have to bound $B_\sigma(x)$, and this has to be done uniformly in $x$. We will split $\R$ into two domains, and show that, in some sense, our choice of $\beta$ already deals with $B_\sigma$ for large $x$, so that the crucial region is near the minimum $x_1$.

\paragraph{What happens for large $x$}%
 \newcommand{\phiPrime}{\Phi'_\sigma(x)}
 \newcommand{\phiSeconde}{\Phi''_\sigma(x)}
 \newcommand{\expMoins}{e^{-\Phi_\sigma(x)}}
 \newcommand{\expPlus}{e^{\Phi_\sigma(x)}}
We study the $x\geq K_2$ by following the proof of corollary 4 in \cite{BCR05}.
\begin{lmm}
  \label{lmm=controleBeta}
For all $\sigma$, there exists a $c_\sigma$ such that:
\begin{equation}%
  \label{equ=cSigma}
  \forall x \notin [K_1,K_2]
  \qquad
  \beta\left(\frac{
    2e^{-\Phi_\sigma(x)}
    }{
    \Phi'_\sigma(x)
    }\right)
    \geq
    \frac{c_\sigma}{\Phi'_\sigma(x)^2}.
\end{equation} %
One may choose $c_\sigma   = \frac{1}{\sigma^2}$.
\end{lmm}
\begin{proof}%
  Recall that the same bound holds for $V$  (\cf{} hypothesis U\ref{hpth=beta}); we try to carry it over to $\Phi_\sigma$.

The behavior of $V$ near its minimum allows us to get an equivalent for $Z_\sigma$ using Laplace's method (\cf{} for example \cite{Die68}); if $V\sim(x-x_1)^b$, we get
\[Z_\sigma \sim C \sigma^{1/b},\]
where $C$ depends only on $V$. Let us bound the argument in the function $\beta$.
\begin{align*}
  2\frac{\exp(-\Phi_\sigma)}{\Phi'_\sigma}
  = \frac{2\sigma}{Z_\sigma} \frac{\exp(-V/\sigma)}{V'}
  \leq C'\sigma^{1-\frac{1}{b}} \frac{\exp(-V/\sigma)}{V'}
  \leq C'\sigma^{1-\frac{1}{b}} \frac{\exp(-V)}{V'}
  \leq \frac{\exp(-V)}{V'}
\end{align*}
for $\sigma$ small enough, because $b$ is strictly greater than $1$, so that $\sigma^{1-1/b}$ goes to zero. 
Since $\beta$ decreases, we get, outside $[K_1,K_2]$:
\[
\beta\left(
  2\frac{\exp(-\Phi_\sigma)}{\Phi'_\sigma}
  \right)
  \geq \beta\left( \exp(-V)/V'\right) \geq \frac{1}{V'^2} 
  = \frac{1}{\sigma^2 \Phi_\sigma'^2}.
  \qedhere
\]
\end{proof} %
\begin{lmm}
  \label{lmm=bornesPourXGrand}
  For all $x\geq K_2$, we have the following inequalities:
  \begin{align*}%
    \int_{m_\sigma}^x e^{\Phi_\sigma} 
    &\leq \int_{m_\sigma}^{K_2} e^{\Phi_\sigma(y)}dy
    + 2\frac{e^{\Phi_\sigma}}{\Phi'_\sigma}\\
    \mu_\sigma([x,\infty))
    &\leq 2\frac{
      e^{-\Phi_\sigma(x)}
      }{
      \Phi'_\sigma(x)
      }
    \leq 3\mu_\sigma([x,\infty)).
  \end{align*}%
\end{lmm}
\begin{proof}
 For all $x\geq K_2$ and $\sigma$ small enough (less than $1/(2C_V)$), the hypothesis on $V$ gives us:
 \begin{align*}%
 \frac{\abs{\Phi''_\sigma(x)}}{\Phi'_\sigma(x)^2} 
   = \sigma \frac{\abs{V''}}{V'^2}
   \leq C_V\sigma \leq  \frac{1}{2}.
 \end{align*}
 Therefore:
 \[ %
 \left(\frac{
   \exp(\Phi_\sigma)
   }{
   \Phi_\sigma'
   }\right)' 
   \geq \frac{1}{2} e^{\Phi_\sigma}, 
 \] %
 This gives the first result by integration. In a similar way, 
 \[%
 \left(\frac{
   \exp( - \Phi_\sigma)
   }{
   \Phi_\sigma'
   }\right)' 
   \in  \left[ \frac{1}{2} e^{ - \Phi_\sigma}, \frac{3}{2} e^{-\Phi_\sigma} \right]
\]%
leads to the second claim.
\end{proof}
We are now in a position to bound $B_\sigma(x)$.
\begin{align*} %
  B_\sigma(x)
  & = \mu_\sigma([x,\infty))\times \frac{1}{\beta(\mu_\sigma([x,\infty)))}
    \times \int_{m_\sigma}^x \expPlus  & \\
  &\leq  2 \frac{\expMoins}{\phiPrime}
    \times \frac{1}{\beta\left( 2 \frac{ e^{-\Phi_\sigma}}{\Phi'_\sigma}\right)}
    \times \left(
    \int_{m_\sigma}^{K_2} e^{\Phi_\sigma(y)}dy
      + 2\frac{\expPlus}{\phiPrime}
      \right) & \text{(by lemma \ref{lmm=bornesPourXGrand})}\\
  &\leq 2 \frac{\expMoins}{\phiPrime}
    \times \frac{\phiPrime^2}{c_\sigma}
    \times \left(
    \int_{m_\sigma}^{K_2} e^{\Phi_\sigma(y)}dy
      + 2\frac{\expPlus}{\phiPrime}
      \right)  & \text{(by lemma \ref{lmm=controleBeta})}\\
    &\leq
    \frac{2}{c_\sigma} \phiPrime 
      \int_{m_\sigma}^{K_2} e^{V_\sigma(y) - V_\sigma(K_2)}dy
      + \frac{4}{c_\sigma} .  & \text{(because $V(x)\geq V(K_2)$)}\\
\end{align*} %
The hypotheses imply that $V(y)\geq V(K_2)$, whenever $\abs{y}\leq K_2$. On the other hand, $\Phi'_\sigma$ is bounded above by $C/\sigma$ (since $V'$ is supposed to be bounded). Finally, 
\[
  \forall x \geq K_2, \qquad B_\sigma(x) \leq \frac{C'}{c_\sigma \sigma}, 
\]
where $C'$ is independent of $\sigma$.

\paragraph{What happens in the well} %
The general strategy here is to bound $B_\sigma(x)$ by studying only the numerator. The denominator can be (very) roughly bounded by $\beta(1/2)$ (which does not depend on $\sigma$). 
The partition function disappears, and we get:
\[
  B_\sigma(x) \leq C \int_{m_\sigma}^x e^{V(y)/\sigma} dy
    \times \int_x^{\infty} e^{-V(y)/\sigma} dy.
\]
We need a bound on $V$ near the median: under our hypotheses, since $\mu_\sigma$ converges weakly to $\delta_{x_1}$, the continuity of $V$ in $x_1$ yields (for $\sigma$ small enough):
\[%
\forall x\in [(x_1,m_\sigma)], \qquad V(x)\leq \dStar/4.
\]
Now we can bound the  first integral in the following way:
\[%
  \int_{m_\sigma}^x e^{V(y)/\sigma} 
  \leq 
    (K_2 - K_1) \exp\left( \frac{1}{\sigma}  \max(s(x),\dStar/4)\right) ,
\]
where $ \dStar/4$ takes care of the case when $m_\sigma$ is less than $x_1$.

We cut the second integral in two parts:
\[
\int_x^\infty e^{-V_\sigma(y)} dy \leq \int_x^{K_2} e^{-V_\sigma(y)}dy + \int_{K_2}^\infty e^{-V_\sigma(y) dy}.
\]
Since $V$ is strictly increasing after $K_2$, we may apply Laplace's method to the second term. In the first one, we use a rough bound on $V$:
\[
\int_x^\infty e^{-V_\sigma(y)} dy
\leq (K_2-K_1) \exp\left(-\frac{i(x)}{\sigma}\right)  
  + C \exp\left(-\frac{V(K_2)}{\sigma}\right).
\]
Since $i(x)$ is less than $V(K_2)$, the second term is less than the first one (up to a constant), and there exists $C'$ such that:
\[
\int_x^\infty e^{-V(y)} dy \leq C' \exp\left( - \frac{i(x)}{\sigma}\right).
\]
Coming back to $B_\sigma$, we get:
\begin{align*}
B_\sigma(x) 
  &\leq C''\exp\left(
    \frac{1}{\sigma}\left( \max(s(x),\dStar/4) - i(x) \right)
  \right)\\
  &\leq C''' \exp\left(\frac{\dStar}{\sigma}\right).
\end{align*}

\paragraph{Conclusion~: An upper bound on $\beta$}
Let us now gather the bounds on $B_\sigma(x)$ we derived in the preceding paragraphs. 
\begin{lmm}
  There exists a $C$ (independent of $\sigma$) such that, for all $\sigma$,
  \[ B_\sigma = \sup_{x\geq m_\sigma} B_\sigma(x)  
  \leq C \exp\left(\frac{\dStar}{\sigma}\right). \]
\end{lmm}
With this result in hand, we may apply Barthe, Cattiaux and Roberto's result (theorem \ref{thm=CapMesBCR}): this  proves the first claim of theorem \ref{thm=alphaUnidim}.

 The modified function $\alpha_\sigma$ is deduced from $\beta_\sigma$ with the help of theorem \ref{thm=deNormeInfaOrlicz} (see below, in section \ref{sec=CapMesEtPoincareFaible}).

  Finally, the growth hypothesis on $\beta$ (U\ref{hpth=croissanceBeta}) guarantees that, near $0$, $\beta$ is bounded by a power of $\ln(1/s)$; this immediately implies the last result, and concludes the proof.

\section{The weak inequality in any dimension}%
\label{sec=multidim}
We now turn to the proof of the weak inequality (the bound in lemma \ref{lmm=majAlpha}) in any dimension. 
We are going to need one more hypothesis on the structure of potential wells, to avoid ``pathological'' cases.

After that, we proceed in several steps. First we recall our aim and explain the main lines of the proof. 
During this proof, a certain ``path'' (in fact, an open set of $\R^d$)  will appear. It will be used to derive a ``capacity-measure'' inequality. Eventually, we will go from this inequality to the one we seek, using a result from next section.

\subsection{The last hypothesis on the potential}%
To write down the last hypothesis we shall make on $V$, we first need a few more notations.

For all $x\in\R^d$, we call $\Gamma_x$ the set of paths from $x$ to $0$. For each such $\gamma$ ($\gamma$ is a continuous function from $[0,1]$ into $\R^d$), we call $h(\gamma)$ the ``height'' of $\gamma$, \ie{} the highest value taken by $V$ along $\gamma$:
\[ h(\gamma) = \sup_{t\in[0,1]} V(\gamma(t)).\]
Now suppose we try to go from $x$ to $0$ while remaining as low as possible (\ie{} we try to find a path where $V$ is small). There is a minimum price to pay; whatever path we choose, we will necessarily go at least as high as:
\[ h(x) = \inf_{\gamma \in \Gamma_x} h(\gamma).\]
We will call ``good paths'' the ones that stay below that minimal height:
\[ \gamma \text{ is good} \Leftrightarrow h(\gamma) = h(x). \]
A priori, for a given $x$, a good path from $x$ to $0$ need not exist: it may well be the case that, if one tries to find $\gamma$ such that $h(\gamma) \leq h(x) + 1/n$, one has to go farther and farther as $n$ grows, and that no finite path achieves the infimum bound. 

Finally, the height of the ``potential barrier'' between $x$ and the global minimum will be called $d^\star(x)$:
\[ d^\star(x) = h(x) - V(x), \]
and the height of the biggest barrier will be just $d^\star$:
\[ d^\star = \sup_{x} d^\star(x).\]

\begin{hpth}%
  \label{hpth=bonsChemins}
  The potential barriers have a bounded height:
  \[ d^\star < \infty.\]
  Moreover, each point can reach $0$ by a  ``relatively short'' good path. More precisely, there exists a function $R$ (a maximal radius), from $\R^d$ to $\R$, which satisfies the following conditions:
  \begin{itemize}
  \item For all $x$, the ball centered in zero and of radius $R(\abs{x})$ contains a good path for $x$:
      \[
	\forall x, \exists \gamma\in\Gamma_x,
	\begin{cases}
	  \gamma([0,1]) &\subset \ball{R(\abs{x})} \\
	  \gamma &\text{is good}.
	\end{cases}
      \]
  \item The function $R$ grows like a power of the distance to the origin:
      \[ R(\abs{x}) \leq c_{R}\abs{x}^{d_R}.\]
  \end{itemize}
\end{hpth}%

\subsection{The one-point weak inequality} %
As was said before, we will not prove in this section a full weak Poincaré inequality, \ie{} we will not get \eqref{eq=defPoincareFaible} for all $r$. Instead, we just prove it for a specific value of $r$, namely the $r =r_t = (\ln t)^{-2} (\ln \ln t)^{-8} $ (\cf{} equation \eqref{eq=choixDeR}).  Since $\sigma(t) = c/(\ln t)$, we note that:
\begin{align*}
  r_t &\geq C \frac{\sigma^2}{\ln (\sigma)^8} \\
  &\geq C' \sigma^m,
\end{align*}
for some $C,C'$ and $\sigma$ small enough. Therefore, and since $\alpha_t$ decreases, we may prove an inequality with $\sigma^m$ instead of $r_t$. 

More precisely, we will get:
\begin{thm}%
  \label{thm=PoincareDefectif}
  Let $m$ be a real number, strictly smaller than $1+\croissanceV$, and let $D^\star$
  be a constant, $D^\star>d^\star$. Then there exists a $C_m$ such that, for all $\sigma$, 
  the measure $\mu_\sigma$ satisfies the following one-point weak Poincaré inequality
  \[
  \forall f, 
  \var_{\mu_\sigma}(f)
  \leq C_m \exp\left( \frac{D^\star}{\sigma} \right) \int \abs{\grad g}^2 d\mu_\sigma
  + \sigma^m \orlPhi{f-m_f}^2.
   \]
   where $m_f$ is a median of $f$ under $\mu_\sigma$.
\end{thm}%
As was noted before, this entails
\[ \alpha_t(r_t) \leq \alpha_t (\sigma^m) \leq C_m \exp\left( D^\star/\sigma \right),\]
which is the result of lemma \ref{lmm=majAlpha}.

\bigskip

The end of the section is devoted to the proof of the theorem. It can be sketched as follows.

The idea is to use a capacity-measure criterion restricted to certain sets (large enough sets). Intuitively, if a set $A$ has a large $\mu_\sigma$ mass, it must contain points near the origin; and these points are the important ones, for measuring capacity as well as mass. For these sets, located near the origin, everything should behave as in the compact case, and the inequality should depend on $\sigma$ in the same way as when a Poincaré inequality holds.

Let us fix $D^\star$, strictly bigger than $d^\star$. As was just said, we would like to compare the capacity and measure of large enough sets: let $\kappa>0$ be the minimum mass we will consider ($\kappa$ will depend on $\sigma$). Let $A$ be a Borel set such that:
\[\mu_\sigma(A) \geq 2\kappa(\sigma).\]
Restricting ourselves to these large sets localizes the problem in some sense. To be more precise, we introduce two radii. 
The first one, $r_\sigma$, is such that:
\[ \mu_\sigma (\ball{r_\sigma}) \geq 1 - \kappa.\]
The second one is deduced from it: it is a radius big enough to include good paths (\cf{} hypothesis \ref{hpth=bonsChemins}) starting from any point in the small ball $\ball{r_\sigma}$.
\[ R_\sigma = R(r_\sigma).\]
These two quantities depend on $\sigma$ and $\kappa$; we will see that, for our choice of $\kappa$, $r_\sigma$ and $R_\sigma$ won't grow too fast as $\sigma$ goes to zero.

Let $A = A'\cup A''$, where $A' = A\cap \ball{r_\sigma}$ and $A''$ is the complement set. Since $\mu_\sigma(A) \geq 2\kappa$ and $\mu_\sigma(A'') \leq \kappa$ (by definition of $r_\sigma$), $\mu_\sigma(A') \geq \kappa$, and:
\[ \mu_\sigma(A)  = \mu_\sigma(A')  + \mu_\sigma (A'') \leq 2\mu_\sigma(A').\]
Intuitively, we need only consider the subset $A'$, because it concentrates enough mass.

At this point, our set $A'$  may still be very complicated. In particular, it could be scattered all over the ball $\ball{r_\sigma}$. To avoid this, we will once again restrict ourselves to a subset, trying to keep enough mass in the process.

This is done by cutting $\ball{r_\sigma}$ into small cubes. The bound on the gradient of $V$ (hypothesis \ref{hpth=gradientBorne}) helps us choose a good mesh, such that $V$ does not vary too much inside a little cube.
\begin{prop}
  For all $\eta$, there exists $\epsilon$ (depending only on $V$ and $\eta$), such that, on each cube $B$ with radius $\epsilon$, 
  \[\sup_B V - \inf_B V \leq \eta.\]
\end{prop}
The parameter $\eta$ will be chosen later. 

So we cut $\ball{r_\sigma}$ into many little cubes of radius $\epsilon$. This requires a certain number of cubes, which we call $n_\sigma$. We then have:
\begin{equation}
  \label{eq=decoupageDesBoules}
  \ball{r_\sigma} = B_1 \cup B_2 \ldots B_{n_\sigma}.
\end{equation}
In the same way, $N_\sigma$ will be the number of cubes necessary to cover $\ball{R_\sigma}$. 
We denote by $A_i$ the intersection of $A$ and $B_i$. We apply the pigeonhole principle to say that one of the $A_i$'s must be large enough:
\[ \exists i_0, \mu_\sigma(A_{i_0}) \geq \frac{1}{n_\sigma} \mu_\sigma(A').\]
To sum up our considerations on sets, for each $A$, we have found a subset $A_{i_0}$ such that:
\begin{itemize}
\item $A_{i_0}$ is a subset of a cube of radius $\epsilon$, 
\item $A_{i_0}$ is not too far from the origin ($A_{i_0}\subset \ball{r_\sigma}$)
  \item $A_{i_0}$ is big enough compared to  $A$~: $\mu_\sigma(A_{i_0}) \geq \frac{1}{2n_\sigma} \mu_\sigma(A)$.
\end{itemize}
In some sense, we need only consider the case when $A$ looks like a ball and is not too far from the origin. 
We are going to see how this can be used to build a certain path between $A_{i_0}$ and $0$, and from this path, deduce a capacity-measure inequality.

\subsection{Building a path and straightening it out}%
\newcommand*{\chemin}{\mathcal{C}_A}
\newcommand*{\droitChemin}{\mathcal{L}_A}
\newcommand*{\constJacobien}{C_\phi}
\newcommand*{\transpose}[1]{{}^t#1}
\newcommand*{\compose}{\circ}
Recall that our goal is to compare the capacity and the measure of sets, and more precisely to bound the capacity from below and the measure from above.

The capacity is defined by an infimum bound:
\begin{equation}
  \label{eq=defCapacite}
  \capacite_\mu(A)  = \inf \left\{ \int \abs{\grad f}^2 d\mu, \ind{A} \leq f \leq 1, \mu(\supp f) \leq \frac{1}{2} 
  \right\} .
\end{equation}
Note that we only define capacities for sets whose measure is less than $1/2$. This restriction explains why we use function recentered by their median when we deduce functional inequalities from capacity-measure criteria.

Since we seek a bound from below, we consider a function satisfying the conditions, and we try to bound:
\[ \int \abs{\grad f}^2 d\mu.\]
The key idea is to find out a region of $\R^d$ which should contribute a lot to this integral. Since the function $f$ equals $1$ near $A$, and $0$ near $0$ (the measure of its support being less than $1/2$), there must be a transition between $A$ and $0$: this is where the gradient of $f$ appears. Still on the intuitive level, if the integral is to be small, we had better make this transition in a region where $\mu$ has less mass, \ie{} in a zone where $V$ is large. This is the reason why we introduced the good paths: to go from $A$ to zero, a large contribution to the energy should appear along these good paths.

To put these ideas on a firm ground, we will build, starting from $A$ (or more precisely from $A_i$), an open set $\chemin$ with good regularity properties, and then bound the capacity by integrals over this open set. This construction is depicted in figure \ref{fig=constructionChemin}.

Once this set is built, we proceed in two steps. First, for all function $f$ satisfying the conditions of \eqref{eq=defCapacite},
\[ \int \abs{\grad f}^2 d\mu_\sigma \geq \int \abs{\grad f }^2 \ind{\chemin} d\mu_\sigma\]
On the path, we know by design that $V$ is bounded above by $V(x_A^\star) + \eta$. Indeed, $V$ is less than $V(x_A^\star)$ along $\gamma$, and the size $\epsilon$ of the cubes has been chosen so that on each cube, the oscillation of $V$ is less than $\eta$. Therefore, we may compare our integral with an integral with respect to the Lebesgue measure.
\begin{equation}
  \label{eq=premiereMinorationCap}
  \int \abs{\grad f}^2 d\mu_\sigma 
  \geq \frac{1}{Z_\sigma}\exp\left( -V(x_A^\star) - \eta\right) \int \abs{\grad f}^2 \ind{\chemin}  d\lambda.
\end{equation}
The next step is to bound the latter integral on $\chemin$. Our only hypotheses is that $f$ must be $1$ on $A$, and $0$ near zero. 
The idea is then to apply a Poincaré inequality to compare the energy to a variance. Unfortunately, though we know that a Poincaré inequality effectively holds under quite general assumptions for a  bounded domain in $\R^d$ (this is proved in many textbooks on partial differential equations, see \eg{} \cite{Eva98},p.~275--276), the explicit constants and their behaviour when the domain changes is not well known. However, there is a case for which we have such explicit estimates, namely the case of convex domains.

\begin{thm}[Poincaré inequality for convex domains]
  \label{thm=poincareConvexe}
  Let $L$ be a convex bounded domain in $\R^d$. Then the Lebesgue measure on $L$ satisfies a Poincaré inequality, and the constant can be bounded above using only the diameter $d_L$ of the domain:
  \[
  \var_{\lambda_L}(f) = \int \left( f - \int fd\lambda_L\right)d\lambda_L
  \leq \frac{d_L^2}{\pi^2} \int \abs{\grad{f}}^2 d\lambda.
  \]
\end{thm}
This theorem is proved \eg{} by Payne and Weinberger, and Bebendorf in \cite{PW60,Beb03}.  Note that other bounds in more complicated cases have been derived (see \cite{CL97} for star-shaped domains, or \cite{Che90} for bounds depending on the geometry of the boundary).
  
  Note that, by abuse of notation,  we use $\var$ for a non-normalized measure.

  \bigskip

In order to use this result, we try to ``straighten out'' the set $\chemin$.

We will build a function $\phi$ sending $\chemin$ to a tube $\droitChemin$. This function will be defined piecewise, on each of the little cubes that $\chemin$ crosses. Let us denote these cubes as $C_0, \ldots C_m$. It is easy to see that the intersection of $\chemin$ and one of these cubes can only take a finite number of shapes (up to a rotation and/or translation). In $d=2$ for example, only two different shapes are possible (either a straight tube or a bended one, see figure \ref{fig=redressement}). Each of these shapes may be ``straightened out'' into a tube by a diffeomorphism. We have to be a bit careful in choosing these diffeomorphisms $\phi_j$ (one for each shape). We will ask two things: they should behave like a rigid motion in the neighborhood of the edges (so we may ``glue'' two transformations together), and their Jacobian matrix should be sufficiently ``nice'' (the ``niceness'' needed will be made precise later). Such a choice is possible; see the figure \ref{fig=redressement} for an explanation of a possible way to find such good functions. 

Once this is done, we only have to glue our pieces together. Let us denote the pieces $\chemin\cap C_i$ by $T_i$. We leave $T_0$ where it stands, and look at $T_1$.  We have seen that it may be straightened into a tube, $T'_1$: define $\phi$ on $T_1$ to be precisely this transformation. Now consider $T_2$: we can straighten it by one of our $\phi_j$, and then use a rotation and/or a translation to put it next to $T'_1$. Since we have asked that the $\phi_j$ should be rigid motions near the edges, the two pieces of $\phi$ define a diffeomorphism from $T_1 \cup T_2$ to the straight tube $T'_1\cup T'_2$. We may iterate the process and eventually we get a diffeomorphism $\phi$ from $\chemin$ to $\droitChemin$. One can see on the figure that a little extra care is needed to deal with the end of the path $\chemin$ --- however, adding just one $\phi_j$ to our set of transformations settles the question.

\begin{figure}
  \caption{Building the path $\droitChemin$}
  \label{fig=constructionChemin}\label{fig=redressement} 
\input{chemins}

\medskip
\input{redressement}
\end{figure}

Remember that our goal is to use the Poincaré inequality on the convex set $\droitChemin$. For this to work, we need to control some quantities related to the map $\phi$.

\begin{prop}%
  \label{prop=controlesDuDiffeo}
  There exists a constant  $\constJacobien$, which may depend on $\epsilon$ but not on $\sigma$, such that, at every point, the Jacobian matrix $J_\phi$ satisfies
  \begin{align*}
    \abs{\det(J_\phi)} &\leq \constJacobien, \\
    \lambda_1(J_\phi \transpose J_\phi) &\geq \constJacobien^{-1},
  \end{align*}
  where $\lambda_1(M)$ is the smallest eigenvalue of the symmetric matrix $M$.
\end{prop}%
\begin{proof} %
  This holds by design of the map $\phi$. At each point, $\phi$ is the composition of a rigid motion (which has no effect on the eigenvalues or the determinant of the Jacobian matrix), and of one of the $\phi_j$. For a given $\phi_j$, the properties hold: we have designed the $\phi_j$ as restrictions of diffeomorphisms on larger sets, so the bounds hold by compacity. Since there is a finite number of $\phi_j$, we may choose bounds that do not depend on $j$. This proves that the bounds hold for $\phi$.
\end{proof}%

We may now give our ``straightening'' its rigorous form, namely a change of variables.
\begin{prop}%
  Let $\mathcal{U}$ and $\mathcal{V}$ be open sets, let $\phi$ be a diffeomorphism from $\mathcal{U}$ onto $\mathcal{V}$. If the inequalities in the preceding lemma hold with a constant $C_\phi$, then for all continuously differentiable function $f$ on $\mathcal{U}$, we have:
  \[ \int_{\mathcal{U}} \abs{ \grad f}^2 d\lambda 
  \geq \constJacobien^{-2} \int_{\mathcal{V}} \abs{\grad g}^2 d\lambda.\]
  where $g = f\compose\phi^{-1}$.
\end{prop}%
\begin{proof}
  It's a change of variables. Let us define $F$ by $F(x) = \abs{\grad{f}}^2(x)$. Then:
  \begin{align*} 
    \int_{\mathcal{U}} \abs{\grad f}^2 dx = \int_{\mathcal{U}} F(x)dx
    &= \int_{\mathcal{V}} F\compose\phi^{-1} \abs{\det J_{\phi}^{-1}} dy \\
    &\geq \frac{1}{\constJacobien} \int_{\mathcal{V}} F\compose\phi^{-1} dy.
  \end{align*}
  Since $f = g\compose\phi$, the gradients are given by:
    \[ (\grad f)_{x} = (\transpose J_\phi)_x (\grad g)_{\phi(x)}.\]
  Taking norms, and using the lower bound on the first eigenvalue, we get:
    \[ (\abs{\grad f}^2)_{x} = \transpose \grad g J_\phi \transpose J_\phi \grad g
     \geq \constJacobien^{-1} (\abs{\grad g}^2)_{\phi(x)}.\]
  Rewriting this in $y$ variables, 
  \[ F\compose \phi^{-1}(y) = (\abs{\grad f}^2)_{\phi^{-1}(y)} \geq \constJacobien^{-1}(\abs{\grad g}^2)_{y}.\]
  Finally, 
  \[ \int_{\mathcal{U}}\abs{\grad f}^2 \geq \frac{1}{\constJacobien^2} \int_{\mathcal{V}} \abs{\grad g}^2.\]
\end{proof}

Putting the last two propositions together, we can show:
\begin{prop}%
  \label{prop=minorationEnergieParMasseDesEnsembles}
  There exists a $C$, depending only on $\epsilon$, such that if $f$ satisfies the following conditions:
  \begin{enumerate}
    \item $f$ is continuously differentiable from $\chemin$ into $[0,1]$, 
    \item $\lambda(\{ f=0\} ) \geq l_0$, 
    \item $\lambda(\{ f = 1 \} ) \geq l_1$, 
  \end{enumerate}
  then
  \[ \int \abs{\grad f}^2 d\lambda \geq \frac{C}{N_\sigma^2} \min(l_0,l_1).\]
\end{prop}%
We recall that $N_\sigma$ is the number of balls of radius $\epsilon$ needed to cover the big ball $\ball{R_\sigma}$.
\begin{proof}%
  Suppose $f$ satisfies the hypotheses. Define $g = f\compose \phi^{-1}$ as in the preceding proposition. The various bounds needed on the Jacobian matrix of $\phi$ are provided by proposition \ref{prop=controlesDuDiffeo}. These bounds also imply that $g$ must vanish at least on a set of Lebesgue measure $C_\phi^{-1} l_0$, the same being true for the set where $g=1$. The change of variables has shown:
  \[ \int \abs{ \grad f}^2 d\lambda \geq C_\phi^{-2} \int \abs{\grad g}^2 d\lambda.\]
  On the right hand side, we can now use the Poincaré inequality:
  \[
    \int \abs{\grad g }^2 d\lambda 
    \geq \frac{1}{C_P(\droitChemin)} \var_\lambda (g).\]
  The very purpose of our change of variables was to make the domain convex, so we could make use of theorem \ref{thm=poincareConvexe}. The constant may therefore be bounded by the square of the diameter of $\droitChemin$. Since $\droitChemin$ results from gluing together at most $N_\sigma$ little cubes of radius $\epsilon$, the square of the diameter may be bounded by $\epsilon^2 N_\sigma^2$.
  
  We now turn to the variance, and use the information on the sets where $g$ is $0$ or $1$. We denote by $l'_0,l'_1$ the respective measures of these sets, and by $m$ the mean of $g$ ($m\in[0,1]$). Then:
  \[
    \var (g)  = \int ( g-m)^2 d\lambda  \geq  m^2 l'_0 + (1-m)^2 l'_1.
  \]
  The right hand side is easily shown to be greater than $(l'_0l'_1/(l'_0 + l'_1)$. The latter is bounded below by half the minimum of $l'_0$ and $l'_1$ (because the numerator is less than $2\max(l'_0,l'_1)$). Since $l'_0 \geq C_\phi^{-1} l_0$, and the similar result holds for $l_1$, 
  \[ \int \abs{\grad f}^2 d\lambda \geq \frac{C_\epsilon}{N_\sigma^2} \min(l_0,l_1).\qedhere \]
\end{proof}%
We may now prove the measure-capacity inequality we are looking for. Indeed, recall that our aim is to bound the capacity of a set $A$ from below by a function of its measure. The previous inequality is almost what we want: on the left hand side is (up to a factor, see \eqref{eq=premiereMinorationCap} above) the quantity whose infimum gives the capacity (equation \eqref{eq=defCapacite}), and on the right hand side $l_0$ and $l_1$ are measures of some sets. It remains to show that these measures may be compared to the measure of $A$.

\subsection{The measure-capacity inequality}%
Let us put together the results from the previous section (equation \eqref{eq=premiereMinorationCap} and proposition \ref{prop=minorationEnergieParMasseDesEnsembles})
\begin{align}
 \notag \int \abs{\grad f}^2 d\mu_\sigma
    &\geq \frac{1}{Z_\sigma}\exp\left( -V(x_A^\star) - \eta\right) \int \abs{\grad f}^2  d\lambda 
    \\
    \label{eq=minorationParMinL0L1}
    &\geq \frac{C_\epsilon}{Z_\sigma N_\sigma^2} \exp\left( -V(x_A^\star) - \eta \right) \min(l_0,l_1),
\end{align}
where $l_0,l_1$ are the Lebesgue measure of the following sets:
\begin{align*}
  l_0 &= \lambda (\{f=0\} \cap \chemin)   &   l_1 &= \lambda (\{f = 1\} \cap \chemin).
\end{align*}

To bound $l_0$, we use the fact that $f$ vanishes on a sufficiently large set (as measured by $\mu_\sigma$). Since $\mu_\sigma$ concentrates around $0$, $f$ should vanish near the origin.  More precisely, for a fixed $\epsilon$, we know that for $\sigma$ small enough, the cube centered in $0$ and of radius $\epsilon$ concentrates $3/4$ of the measure. If this cube is labelled $B_0$, we have:
\[ \mu_\sigma(\{f = 0 \} \cap B_0) \geq \frac{1}{4}.\]
Since $V$ is non negative, $\mu_\sigma$ and $\lambda$ are easily compared.
\begin{align*}
\frac{1}{4} \leq \mu_\sigma(\{f = 0\} \cap B_0) 
    &=    \frac{1}{Z_\sigma} \int \ind{f=0}\ind{B_0} \exp( - \frac{V}{\sigma} ) d\lambda \\
    &\leq \frac{1}{Z_\sigma} \int \ind{f=0}\ind{B_0} d\lambda
\end{align*}
The integral on the right hand side is less than $l_0$, 
therefore:
\[ l_0 \geq m_0 =  \frac{Z_\sigma}{4}.\]

Let us derive a similar bound, $m_1$,  for $l_1$. On the cube $B_{i_0}$, $V\geq V(x_A) - \eta$, so:
\begin{align*}
  \mu_\sigma(A_{i_0}) &= \frac{1}{Z_\sigma}\int \ind{A_{i_0}} \exp (-\frac{V}{\sigma}) d\lambda \\
  &\leq \frac{1}{Z_\sigma} \int \ind{A_{i_0}} \exp( -\frac{V(x_A)}{\sigma} + \frac{\eta}{\sigma} ) d\lambda \\
  &\leq \frac{1}{Z_\sigma} \exp \left( \frac{-V(x_A) + \eta}{\sigma}\right) \lambda(A_{i_0}).
\end{align*}
Therefore:
\begin{equation}
  \label{eq=minorationL1}
  l_1 \geq \lambda(A_{i_0}) 
  \geq m_1 = Z_\sigma \exp\left( \frac{V(x_A)}{\sigma} - \frac{\eta}{\sigma} \right) \mu_\sigma(A_{i_0}).
\end{equation}

Since we would like to control $\min(l_0,l_1)$, we now have to compare the two bounds $m_0$ and $m_1$. This is possible thanks to the following inequality:
\[ \mu_\sigma(A_{i_0}) \leq \frac{1}{Z_\sigma} \exp \left( \frac{ - V(x_A) + \eta }{\sigma} \right) \epsilon^d.\]
If we gather almost all terms on the left hand side, we recognize $m_1$: %
\[ m_1 \leq \epsilon^d.\]
Since $m_0 = Z_\sigma/4$, it holds that $m_0 \geq Z_\sigma m_1 \epsilon^{-d}$, and since $Z_\sigma$ goes to zero, it also holds that $m_1 \geq Z_\sigma m_1 \epsilon^{-d}$, so that both $l_0$ and $l_1$ may be bounded below by this quantity:
\[ \min(l_0,l_1) \geq \frac{Z_\sigma^2}{\epsilon^d} \exp\left(\frac{V(x_A)}{\sigma} - \frac{\eta}{\sigma} \right) \mu_\sigma(A_{i_0}).\]
Going back to \eqref{eq=minorationParMinL0L1}, we conclude:
\begin{align*}
  \int \abs{\grad f}^2 d\mu_\sigma 
    &\geq \frac{C_\epsilon}{Z_\sigma N_\sigma^2} \exp\left( -V(x_A^\star) - \eta \right) \min(l_0,l_1) \\
    &\geq \frac{C'_\epsilon Z_\sigma}{N_\sigma^2} \exp\left(V(x_A) - V(x_A^\star) -2\eta \right) \mu_\sigma(A_{i_0}).
\end{align*}
By definition of $x_A^\star$,  $V(x_A) - V(x_A^\star) - 2\eta \geq - d^\star -2\eta \geq -D^\star$. 
On the other hand, $A_{i_0}$ was chosen precisely because it contained enough of $A$'s mass:
$\mu_\sigma(A_{i_0}) \geq (2n_\sigma)^{-1}\mu_\sigma(A)$.
Finally, every function $f$ we can choose in the definition of capacity must satisfy:
\[
  \int \abs{\grad f}^2 d\mu_\sigma 
  \geq \frac{C'_\epsilon \sigma^2}{N_\sigma^2 n_\sigma} \exp\left( -\frac{D^\star}{\sigma}\right) \mu_\sigma(A),
\]
where $C'_\epsilon = CC_\epsilon$.
Taking the infimum over all possible $f$ finally yields the following result.
\begin{prop}
  Let $\kappa(\sigma)$ be a positive number, less than $1/2$. Let $n_\sigma,N_\sigma$ be defined as in the discussion near equation \eqref{eq=decoupageDesBoules}. Then the following bound holds:
\begin{equation}
  \label{eq=capMesIntermediaire}
  \forall A, \mu_\sigma(A) \geq \kappa(\sigma)
  \quad \implies \quad
  \mu_\sigma(A) \leq
    \frac{N_\sigma^2 n_\sigma}{C'_\epsilon\sigma^2}  
    \exp\left( \frac{D^\star}{\sigma} \right) \capacite_{\mu_\sigma}(A).
\end{equation}
\end{prop}

\subsection{Conclusion}%
The bigger part of the proof has now been done; the last thing we need to check is that the number of balls $n_\sigma$ and $N_\sigma$ do not grow too fast as $\sigma$ decreases. Then we will apply theorem \ref{thm=deCapMesAPoincareDefectif} to deduce the one-point inequality of theorem \ref{thm=PoincareDefectif} from our measure-capacity inequality.

Recall that we are given a real number $m$, strictly smaller than $1+\croissanceV$.
Define $\kappa(\sigma) = \exp\left( - \frac{1}{\sigma^m} \right)$. We want to find an $r_\sigma$ such that the mass of $\ball{r_\sigma}$ is greater than $1-\kappa$. For any set $A$, we may write:
\begin{align*}%
\mu_\sigma(A)  
&= \frac{1}{Z_\sigma} \int \ind{A} \exp\left( -\frac{V}{\sigma} \right) \\
&= \frac{Z_{2\sigma}}{Z_\sigma} 
  \times \frac{1}{Z_{2\sigma}} \int \ind{A} \exp\left( -\frac{V}{2\sigma} -\frac{V}{2\sigma} \right) \\
&= \frac{Z_{2\sigma}}{Z_\sigma} 
  \times  \int \ind{A} \exp\left( -\frac{V}{2\sigma} \right) d\mu_{2\sigma}.
\end{align*}%
If $V$ takes large values on $A$, we can get a good bound:
\[
\mu_\sigma(A) \leq \frac{Z_{2\sigma}}{Z_\sigma} \exp\left( - \frac{\inf_A V}{2 \sigma} \right) \mu_{2\sigma}(A).
\]
We get rid of the $\mu_{2\sigma}(A)$ by roughly bounding it by $1$. Then we use the growth hypothesis on $V$ (\ref{hpth=croissanceDeV}), with $A = \compl{\ball{r_\sigma}}$. In this case:
\[ \inf_A V \geq \ln(r_\sigma)^{\croissanceV }.\]
We fix an $m'\in ]m,1+\croissanceV[$, and choose:
\[ r_\sigma = \exp\left( \left(\frac{1}{\sigma}\right)^{(m'-1)/\croissanceV } \right), \]
which ensures:
\begin{align*}
  \inf_A V
    &\geq \left( \frac{1}{\sigma}\right)^{(m'-1)}, \\
  \mu_\sigma( \compl{\ball{r_\sigma}} )
    &\leq \frac{Z_{2\sigma}}{Z_\sigma} \exp\left( -\frac{1}{2\sigma^{m'}} \right)  
\end{align*}
The asymptotic behavior of $Z_\sigma$ (\cf{} annex \ref{sec=fonctionDePartition}) implies that $Z_{2\sigma}/Z_\sigma$ converges, and since $m'>m$,
\[ 
\mu_\sigma( \compl{\ball{r_\sigma}} ) \leq \exp\left( - \frac{1}{\sigma^m} \right) 
\]
for $\sigma$ small enough. This shows that $r_\sigma$ satisfies the condition we wanted.

We may now end the proof of the theorem. Coming back to the measure-capacity inequality \eqref{eq=capMesIntermediaire}, 
we note that $R_\sigma$, $n_\sigma$ and $N_\sigma$ all behave like $r_\sigma$ to a certain power
(for $R_\sigma$ we use hypothesis \ref{hpth=bonsChemins}, and $n_\sigma, N_\sigma$ are just a
 number of cubes of fixed radius in the big cubes of side length $r_\sigma$ and $R_\sigma$).
 Therefore, there exists a $C$ such that
\begin{equation}%
  \forall A, \mu_\sigma(A) \geq \kappa(\sigma)
  \quad \implies \quad
  \mu_\sigma(A) \leq
  \frac{r_\sigma^C}{\sigma^2}  \exp\left(
    \frac{D^\star}{\sigma}
  \right) \capacite_{\mu_\sigma}(A).
\end{equation}%
The value of $r_\sigma$ and the fact that $m'-1$ is strictly less than $\croissanceV$ makes $\exp(D^\star/\sigma)$ the 
biggest term, so that, up to a slight increase of $D^\star$, 
\[
  \forall A, \mu_\sigma(A) \geq \kappa(\sigma)
  \quad \implies \quad
  \mu_\sigma(A) \leq
  \exp\left(\frac{D^\star}{\sigma}\right)\capacite_{\mu_\sigma}(A).
\]
This inequality, thanks to theorem \ref{thm=deCapMesAPoincareDefectif} below, implies precisely the one-point weak Poincaré inequality we claimed in theorem \ref{thm=PoincareDefectif}.

\section{A measure-capacity criterion for one-point weak Poincaré inequalities} %
\label{sec=CapMesEtPoincareFaible}
\subsection{Definitions}%
In this section we study the interplay between weak Poincaré inequalities  and measure-capacity inequalities. 
Let us start by recalling exactly what a weak Poincaré inequality is.
\begin{dfn} [M. Röckner and F.Y.~Wang, \cite{RW01}]
  Let $\mu$ be a measure and $\rndNrm$ be a norm, stronger than the $L^2(\mu)$ norm. %
  The measure $\mu$ is said to satisfy a \term{weak Poincaré inequality} for the norm $\rndNrm$
  if there exists a decreasing positive function $\alpha$, defined on $\R_+^\star$ such that:
  \[\forall f\in L^2(\mu), f\text{ such that $\mu f = 0$},
  \forall r>0, \mu(f^2) \leq \alpha(r)\nrj{f,f} + r\rndNrm(f)^2\]
  If this holds, $\alpha$ will be called a \term{compensating function}.
\end{dfn}
\begin{rmq}[on means and medians]
  The original statement on weak Poincaré inequalities involves functions recentred by their mean value $\mu(f)$, and an $L^\infty$ norm. However, the approach by measure-capacity inequalities developed in \cite{BCR05,BCR05b} works with functions recentred by their median $m_f$. When the norm is the sup norm, it is easy to go from one to the other: the three quantities $\osc(f),\ninf{f - m_f}$ and $\ninf{f-\mu(f)}$ are within (universal) bounds of each other. 

  Since we need to work with another norm, we will show that we can still go from $\rndNrm(f - m_f)$ to $\rndNrm(f-\mu f)$ (\cf{} equation \eqref{eq=medianeEtMoyenne} in annex \ref{sec=OrliczNorms}).
\end{rmq}
This is equivalent to the slightly modified definition:
\begin{prop}
  A weak Poincaré inequality holds if and only if:
\begin{equation}\label{eq=WPI}
\forall r>0, \exists c_r, \forall f \in L^2(\mu), \qquad \mu(f) = 0
\implies \mu(f^2)\leq c_r \Nrj{\mu}{f,f}+ r \rndNrm(f)^2.
\end{equation}
If the inequality holds for a given couple $(r,c_r)$, we will say that $\mu$ satisfies a \term{one-point weak Poincaré inequality}.
\end{prop}
Therefore the weak Poincaré inequality holds if and only if a one-point inequality holds for each point $r$.
\begin{proof}
The only thing to check is that we can deduce the inequality of the definition 
from \eqref{eq=WPI}. To each $r$, we associate $c_r$ according to \eqref{eq=WPI}. 
Then we just define $\alpha(r)= \inf\{ c_s ; s\leq r \}$. The function $\alpha$ is decreasing.
Now let $f$ be a function in $L^2$ and $r>0$. For any $\epsilon$, we may find an $s\leq r$ such that:
\[c_s \leq \alpha(r) + \epsilon.\]
If we apply \eqref{eq=WPI} with this $s$, we get (since $s\leq r$):
\begin{align*}
\mu(f^2) &\leq c_s \nrj{f} + s \rndNrm(f)^2\\
&\leq \alpha(r)\nrj{f} + r \rndNrm(f)^2 + \epsilon \nrj{f}.
\end{align*}
Since this is true for any $\epsilon$, we may let it go to zero, and we have found a function $\alpha$.
\end{proof}

We will be specifically interested in these inequalities for one special norm. We now define this norm and recall some of its properties, without proofs. For a short introduction (with the results we need here), see \eg{} \cite{Ale04}; for an extensive treatment we refer to \cite{RR91}.

  Let $\phi,\psi$ be defined on $\R_+$ by $\psi(x) = x \log(1+x)$, $\phi(x) = \psi(x^2)$.
  For any measurable $f$, define the \term{Orlicz norm} (usually called the \term{Luxembourg norm}; there is another natural norm on the Orlicz space, but we won't need it here) of $f$ to be:
  \[
  \orlPhi{f} = \inf\left\{ \lambda, \int \phi\left( \frac{\abs{f}}{\lambda}\right) \leq 1 \right\}.
  \]
  Note that, with this definition, $\orlPhi{1} \neq 1$.  The set of functions $f$ for which this norm is finite is denoted $L_\phi$, it is a vector space, 
  and it is complete for the Orlicz norm.
  In the same way, if $\psiStar,\phiStar$ are the convex dual functions of $\psi,\phi$, we may define the corresponding Orlicz spaces. It is easily seen that for every positive $f$, $\orlPsi{f^2} = \orlPhi{f}$. The dual functions allow us to state the following Hölder-like property.
  \begin{prop}[Hölder-Orlicz]
    If $f,g$ are two measurable functions, respectively in $L_\psi$ and $L_\psiStar$, then $fg$ is in $L^1$, and
    \[ \left| \int fg d\mu \right| \leq 2 \orlPsi{f}\orlPsiStar{g}. \]
  \end{prop}
  The constant $2$ is necessary because we work with Luxembourg norms.
  To conclude this account on Orlicz norm, we recall here the norm of an indicator function:
  \begin{prop}
    Let $A$ be a measurable set. Then $\ind{A}$ is in the Orlicz space $L_\psi$ and:
    \[ \orlPsiStar{\ind{A}} = \hat{\psi}(\mu(A)),\]
    where $\hat{\psi}(x) = \frac{1}{(\psi^\star)^{-1}(1/x)}$.
    Moreover, for all $x$ sufficiently small, we have  the following bound:
    \[ \hat{\psi}(x)\leq \frac{2}{\log(1/x)}.\]
  \end{prop}
  \begin{proof}
    Once again we refer to \cite{Ale04,RR91} for the first result. The explicit bound on $\hat{\psi}$ follows easily from the bound $\psi^\star \leq xe^x$ and the definition of $\hat{\psi}$.
  \end{proof}

\subsection{Measure-capacity inequalities for large sets and one-point inequalities }%
Here we show the result which was used in the preceding section: if we can compare the measure and the capacity of large sets, we can deduce a one-point weak inequality.
\begin{thm}%
  \label{thm=deCapMesAPoincareDefectif}
  Suppose that there exists $\kappa<1/2$, and a real constant $C_\kappa$ such that,  for every set $A$ whose measure is larger than $\kappa$, we have:
  \begin{equation}
    \label{eq=CapMesPourGrandsEnsembles}
    \capacite_\mu(A) \geq C_\kappa \mu(A).
  \end{equation}
  Then $\mu$ satisfies the one-point weak Poincaré inequality:
  \[ \var_\mu(g) \leq \frac{c}{C_\kappa} \int \abs{\grad g}^2 d\mu + \kappa  \osc^2(g),\]
  where $c$ is universal. We may replace the $L^\infty$ norm by an Orlicz norm, in which case the inequality reads:
  \[ \var_\mu(g) \leq
    \frac{c}{C_\kappa} \int \abs{\grad g}^2 d\mu
    + \hat{\psi}(\kappa)  \orlPhi{g-m_g}^2.
   \]
\end{thm}%
\begin{rmq}
  Note that if  \eqref{eq=CapMesPourGrandsEnsembles} holds for \emph{all} sets, regardless of their measure, then $\mu$ satisfies a (strong) Poincaré inequality (since we may take $\kappa = 0$). This is well-known, \cf{} \cite{BCR05b} and references therein. This characterization of a functional inequality in terms of a relation between measures and capacities of sets is in fact more general, and provides a way to compare many functional inequalities. For a detailed account on these questions, and links with isoperimetric properties, we refer to \cite{BCR05b} (especially section 5).
\end{rmq}

\begin{proof}%
  We follow the proof of theorem 2 in \cite{BCR05} (which deals with the (full) weak inequality).

  Let $f$ be a function and $m$ a median for $f$. We cut the space in half, according to whether $f$ is greater than $m$ or not; we denote by $\Omega_+, \Omega_-$ the two sets. The integral may be written as:
  \[ \var_\mu(f) \leq \int(f-m)^2 d\mu  = \int_{\Omega_+} (f-m)^2 d\mu + \int_{\Omega_-} (f-m)^2 d\mu.\]
  We will show how to deal with the leftmost integral, the other one being similar.
  \[c = \inf\{t\geq 0, \mu(g^2 >t) < \kappa\}. \]
  If $c$ is zero, then $\mu(g>0)$ is less than $\kappa$, and:
  \[%
  \int_{\Omega_+}g^2 d\mu \leq 
  \begin{cases} \kappa \sup g^2 \text{ in the $L^\infty$ case,}\\
    \hat{\psi}(\kappa)  \orlPhi{f-m}^2 \text{ in the Orlicz case,}
  \end{cases}
  \]%
  so the inequalities we are looking for hold in the half-space $\Omega_+$.

  Thus we need only consider the case where $c$ is strictly positive. By a continuity argument ($\mu$ will always have a density), we can find a set $\Omega_0$ such that $\mu(\Omega_0) = \kappa$ and $\{g^2>c\} \subset \Omega_0 \subset \{g^2 \geq c\}$. We fix a $\rho>1$, and introduce the level sets $\Omega_k = \{ g^2 \geq \frac{c}{\rho^k}\}$.
  We decompose the integral over these sets:
  \begin{align*}%
  \int_{\Omega_+} g^2  
  &=    \int_{\Omega_0} g^2 d\mu + \sum_{k\geq 1} \int_{\Omega_k \setminus \Omega_{k-1}} g^2 d\mu \\
  &\leq \int_{\Omega_0} g^2 d\mu 
    + \sum_{k\geq 1} \frac{c}{\rho^{k-1}}\left( \mu(\Omega_k) - \mu(\Omega_{k-1}) \right)
  \end{align*}%
  The sum is dealt with thanks to an Abel transform:
  \begin{align*}%
    \sum_{k\geq 1} \frac{1}{\rho^{k-1}}(\mu_k - \mu_{k-1}) \\
  &= \sum_{k\geq 1} \frac{\mu_k}{\rho^{k-1}} - \sum_{k\geq 0} \frac{\mu_k}{\rho^k} \\
  &= \sum_{k\geq 1} \mu_k \left( \frac{1}{\rho^{k-1}} - \frac{1}{\rho^k} \right) - \mu_0.
\end{align*}%
  This is where we do not follow \cite{BCR05}: since we simply suppose an inequality between capacity and measure, 
  we can get rid of the $\mu_0$ and write
  \[ %
  \sum_{k\geq 1} \frac{1}{\rho^{k-1}}(\mu_k - \mu_{k-1}) 
  \leq (\rho-1) \sum_{k\geq 1} \frac{\mu_k}{\rho^k}.
  \]%
  The rest of the proof follows the same line as in \cite{BCR05} --- at this point, we use the measure-capacity inequality on each set $\Omega_k$. They are designed to have their measure bigger than $\kappa$, so that we may apply our hypothesis:
  \[\mu_k \leq \frac{1}{C_\kappa} \capacite(\Omega_k).\]
  Now, to bound the capacity from above, we apply the definition with well-chosen functions $g_k$:
  \[%
  g_k = \min\left(
    1,
    \left( \frac{ g - \sqrt{ c\rho^{-k-1}}}{\sqrt{c\rho^{-k}} - \sqrt{c\rho^{-k-1}}}\right)_+
  \right)
  \]%
  This entails:
  \begin{align*}%
  \mu_k &\leq \frac{1}{C_\kappa} \int \abs{\grad g_k}^2 d\mu  \\
    &\leq \frac{\rho^{k+1}}{C_\kappa c(\sqrt{\rho} - 1)^2 } \int_{\Omega_{k} \setminus \Omega_{k-1}} \abs{\grad g}^2 d\mu.
  \end{align*}%
  Summing over $k$, we get:
  \[%
  \int_{\Omega_+} g^2 d\mu
  \leq \int_{\Omega_0} g^2 d\mu + \frac{\rho(\rho -1)}{ C_\kappa  (\sqrt{\rho} -1)^2} \int \abs{\grad{g}}^2 d\mu.
  \]%
  We may now choose $\rho$; the (non optimal) choice $\rho=4$ gives:
  \[
  \int_{\Omega_+} g^2 d\mu
  \leq \int_{\Omega_0} g^2 d\mu 
    + \frac{12}{ C_\kappa } \int \abs{\grad{g}}^2 d\mu.
  \]
  The only thing left to do is to take care of the integral on $\Omega_0$. This is done with an Hölder-like inequality.
  In the Orlicz norm case, for example, we write:
  \begin{align*}
    \int_{\Omega_0} g^2 d\mu &\leq 2 \orlPsi{g^2} \orlPsiStar{\ind{\Omega_0}} \\
  &\leq 2 \orlPsi{ (f-m)^2 } \hat{\psi}(\kappa) \leq 2\hat{\psi}(\kappa) \orlPhi{f-m}^2.
\end{align*}
  thanks to the Hölder-Orlicz inequality and the relation between $\phi$ and $\psi$ (see the beginning of this section).
\end{proof}%

\subsection{Weak inequalities for different norms} %
To conclude this section, let us state a corollary to the previous result, and prove that weak Poincaré inequalities for many different norms are in fact equivalent. Moreover, if a compensating function is known for one norm, we can immediately deduce a function for another norm; this result was used in the one dimensional case (section \ref{sec=unidim}) where the explicit Hardy-like criteria were known for the $L^\infty$ norm.
\begin{thm}
  \label{thm=deNormeInfaOrlicz}
  Let $\phi,\psi$ be two Young functions, with $\phi(x) = \psi(x^2)$. A measure $\mu$ satisfies a weak Poincaré inequality with the $L^\infty$ norm if and only if it satisfies one with the Orlicz norm $\orlPhi{\cdot}$.

  Moreover, if $\beta$ is a compensating function for the $L^\infty$ norm, then the following function may be chosen for the Orlicz norm:
  \[ \alpha(s)  =  \frac{c}{4} \beta\left(\frac{1}{4} \hat{\psi}^{-1}\left( \frac{s}{2} \right) \right),\]
  where $c$ is universal (and the same as in the preceding result).
\end{thm}
\begin{proof}
  \newcommand{\capMes}{\text{\textbf{M-C}}}
  \newcommand{\poincDef}{\text{\textbf{PWP}}}
  \newcommand{\poincFaible}{\text{\textbf{WP}}}
  First, let us introduce a few notations. We will denote by $\capMes(\kappa,C(\kappa))$ the following comparison between measure and capacity:
  \[ \forall A, \mu(A)>\kappa \implies \capacite(A) \geq C(\kappa) \mu(A).\]
  Similarly,  $\poincDef(r,C(r),\rndNrm)$ will denote the one-point weak Poincaré inequality for a norm $\rndNrm$ with constants $(r,C(r)$, and $\poincFaible(\alpha,\rndNrm)$ will be the (full) weak inequality, with a norm $\rndNrm$ and a compensating function $\alpha$.
  In the previous section, we showed:
  \begin{align*}
  \capMes(\kappa,C(\kappa)) \implies \poincDef\left(\kappa,\frac{c}{C_\kappa},\ninf{\cdot}\right), \\
  \capMes(\kappa,C(\kappa)) \implies \poincDef\left(2\hat{\psi}(\kappa),\frac{c}{C_\kappa},\orlPhi{\cdot}\right).
  \end{align*}
  Going the other way around is easy. Indeed, suppose that $\poincDef(r,C(r),\ninf{\cdot})$ holds. Let $A$ be a set whose measure is less than $1/2$, but greater then $4r$. Let $g$ be any function which may appear in the definition of the capacity of $A$ (\cf{} \eqref{eq=defCapacite}), and let $m_g$ be a median of $g$. Then:
  \[
  \var_\mu g\leq  C_r \int \abs{\grad g}^2 d\mu  + r \ninf{g - m_g}.
  \]
  Without loss of generality, we suppose that $0\leq g\leq 1$, so that the $L^\infty$ norm is bounded by $1$. Moreover, $r$ is less than $\mu(A)/4$, and the variance on the left hand side is bounded below by $(1/2)\min(\mu(A),1/2) \geq (\mu(A)/2)$ (by the same argument used previously, during the proof of proposition \ref{prop=minorationEnergieParMasseDesEnsembles}). This entails:
  \[ \frac{\mu(A)}{2} \leq C_r \int\abs{\grad g}^2 d\mu  + \frac{\mu(A)}{4}.\]
  This immediately implies the measure capacity inequality $\capMes(4r,4/C_r)$.

  If we now try to derive an inequality with an Orlicz norm starting from one with an $L^\infty$ norm, we just translate them in terms of measure and capacity:
  \begin{align*}
    \poincDef(r,C_r,\ninf{\cdot}) &\implies \capMes(4r,4/C_r) \\
    & \implies \poincDef( 2\hat{\psi}( 4r), \frac{c C_r}{4} ).
  \end{align*}
  If we are looking for a full weak Poincaré inequality, we fix an $s$, and define $r=(1/4)\hat{\psi}^{-1}(s/2)$. We may then apply 
  $\poincDef(r,\beta(r))$ to obtain:
  \[ \poincDef(s, c\beta(r)/4,\ninf{\cdot}).\]
  Since $s$ is arbitrary, this concludes the proof.
\end{proof}

\appendix

\section{Orlicz norms, entropy and centering}%
\label{sec=OrliczNorms}
The proof of weak Poincaré inequalities starting from measure-capacity comparisons for an Orlicz norm leads us to consider norms of functions recentered by their median. In fact, what one obtains when applying these criteria is of the form:
\[\var_{\mu}{f^2} \leq \beta(s) \nrj{f} + s \orlPhi{f - m_f}^2,\]
where $m_f$ is a median for $f$. The aim of this section is to bound this term by more tractable quantities (we will use an entropy and a moment).

More precisely we prove the following result:
\begin{lmm}
  \label{lmm=controleOrliczEntropie}
  There exists a $C$ such that, for any positive $f$ and any probability measure $\mu$, the
  following holds:
  \[\orlPhi{f-m_f}^2 \leq C\left( \Ent_\mu(f^2) + 3 \esp_\mu(f^2)\right).\]
\end{lmm}

The proof is done in several steps, and borrows several arguments from \cite{BG99}. First of all, we get rid of the median and replace it by a mean value.
\begin{align}
  \orlPhi{f-m_f} &\leq \orlPhi{f-\mu f} + \orlPhi{\mu f - m_f}\notag\\
  &\leq \orlPhi{f-\mu f} + \abs{\mu f - m_f}.\label{equ=fMoinsMediane}
\end{align}
Let us consider the last term.
\[%
  \mu f - m_f = \int f(x)d\mu - m_f = \int( f-m_f)_+ d\mu - \int(f-m_f)_- d\mu, 
\]
where the integrals are both positive. The absolute value of the left hand side may then be bounded above:
\[
\abs{\mu f - m_f} \leq \max\left( \int(f-m_f)_+ d\mu, \int(f-m_f)_- d\mu \right)
\]
Each of the arguments in the $\max$ can be controlled by Hölder's inequality.
\begin{align*}%
\int (f-m_f)_+ d\mu & = \int( f-m_f) \ind{f>m_f} d\mu \leq \norme{2}{f-m_f} \norme{2}{ \ind{f>m_f} }\\
&\leq \frac{1}{\sqrt{2}} \norme{2}{f-m_f}   & \text{(since  $\mu(f>m_f) < 1/2$)}\\
&\leq \frac{1}{\sqrt{2}}\frac{\sqrt{5}}{2} \orlPhi{f-m_f} &\text{(\cf{} \cite{BG99}, lemma 4.3)}
\end{align*}%
Coming back to \eqref{equ=fMoinsMediane}, we get:
\[\orlPhi{f-m_f} \leq \orlPhi{f-\mu_f} + \abs{\mu f - m_f}
\leq \orlPhi{f-\mu f} + \sqrt{\frac{5}{8}} \orlPhi{f-m_f}.
\]
Since $\sqrt{\frac{5}{8}}\leq 1$, we may put it on the other side to get:
\begin{equation}
  \label{eq=medianeEtMoyenne}
  \orlPhi{f-m_f} \leq C \orlPhi{f-\mu f}
\end{equation}
where $C = (1-\sqrt{\frac{5}{8}})^{-1}$ is universal.

The next step is to bound the Orlicz norm by an entropy. Once again,  we use a result from Bobkov and Götze (\cite{BG99}):
\[ \orlPhi{f-\mu f}^2 \leq \frac{3}{2} \sup_{a\in\R}{\Ent_\mu((f+a)^2)}. \]

Since we would like to deal only with the entropy of $f^2$, we try to compare the entropies of translated functions. Rothaus' lemma tells us:
\[
  \Ent_\mu((f+a)^2) \leq \Ent_\mu(\tilde{f}^2) + 2 \var_\mu(f), 
\]
where $\tilde{f}$ is the centered function $f - \mu f$. The only thing left to do is to bound the entropy of the square of this centered function. This is done in the following lemma.
\begin{lmm}
  \label{lmm=entropieRecentree}
  Let $f$ be a positive function, and $\tilde{f} = f- \mu f$. Then the following holds:
  \[ \Ent_\mu(\tilde{f}^2) \leq \Ent_\mu(f^2) + \int f^2d\mu.\]
\end{lmm}
\begin{proof}%
  Both sides of the equation are homogeneous (of order two), so we may as well suppose $\int f^2 d\mu = 1$. We rewrite the left hand side.
\begin{align*}%
  \Ent_\mu(\tilde{f}^2)  
  &= \int \tilde{f}^2 \log(\tilde{f}^2) d\mu - \esp_\mu(\tilde{f}^2) \log(\esp_\mu(\tilde{f}^2))\\
  &= \int \tilde{f}^2 \log(\tilde{f}^2) d\mu - \var_\mu(f) \log(\var_\mu(f)).
\end{align*}%
The second term is easily dealt with. Indeed, since $\int f^2 =1$, $\var_\mu f$ must be between $0$ and $1$. Since 
$x\mapsto \abs{x\log(x)}$ is bounded by $1/e$ on this interval, one can write:
\[\Ent_\mu(\tilde{f}^2) \leq \int(\tilde{f}^2\log(\tilde{f}^2) d\mu + \frac{1}{e}.\]
   We decompose the integral in two parts, according to whether $f$ is less than $1$ or not.
 \begin{align*}
   \Ent_\mu(\tilde{f}^2) &\leq 
   \int\tilde{f}^2\log(\tilde{f}^2) \ind{\abs{\tilde{f}}\leq 1} d\mu
   + \int\tilde{f}^2\log(\tilde{f}^2) \ind{\abs{\tilde{f}} >}1d\mu + \frac{1}{e} \\
 &\leq 
 \int\tilde{f}^2\log(\tilde{f}^2) \ind{\abs{\tilde{f}} >1} d\mu + \frac{1}{e}, 
 \end{align*}
 since the first term is less than $0$. Now, on the set where $\abs{\tilde{f}}$ exceeds one, $f$ must
 be above its mean: $f$ is indeed positive, and since $\int f^2 d\mu = 1$, $\mu f$ must be in $[0,1]$.
 So $\abs{f-\mu f}$ may be greater than $1$ only when $f$ itself is greater than $1$. 
 This shows that, on $\{\abs{\tilde{f}}>1\}$,
 \[ 1\leq \tilde{f} =  f-\mu f \leq f.\]
 Since $x\mapsto x\log(x)$ increases on $[1,\infty)$, we have:
 \begin{align*}
 \Ent_\mu(\tilde{f}^2) &\leq  \int \tilde{f}^2 \log(\tilde{f}^2) \ind{\abs{\tilde{f}>1}}d \mu + \frac{1}{e}\\
 &\leq \int f^2 \log(f^2) \ind{\abs{\tilde{f}>1}} d\mu  + \frac{1}{e}.
 \end{align*}
 At this point, remark that on $\{ f>1 \}$, $f^2\log(f^2)$ is positive, and since
 $\ind{\abs{\tilde{f}}>1} \leq \ind{f>1}$, 
 \begin{align*}%
 \Ent_\mu(\tilde{f}^2)
 &\leq \int f^2 \log(f^2) \ind{f>1} d\mu + \frac{1}{e}\\
    &\leq \Ent(f^2) - \int f^2\log(f^2)\ind{f<1} d\mu +\frac{1}{e}\\
    &\leq \Ent(f^2) + \frac{2}{e}.
    \end{align*}%
    Since  $\frac{2}{e}\leq 1$, the proof is complete.
\end{proof}%
Gathering our results, we have shown that:
\begin{align*}%
  \orlPhi{f- m_f}^2
  &\leq C\orlPhi{f-\mu_f}^2 & \text{(inequality \eqref{eq=medianeEtMoyenne})}\\
  &\leq \frac{3C}{2} \sup_{a\in\R} \Ent((f+a)^2)  & \text{(Bobkov and Götze's lemma) }\\
  &\leq \frac{3C}{2} \left(\Ent(\tilde{f}^2) + 2 \var_\mu(f)\right)   & \text{(Rothaus' lemma)}\\
  &\leq \frac{3C}{2} \left(\Ent(f^2) + 3\esp_\mu(f^2)\right).   & \text{(lemma \ref{lmm=entropieRecentree})}
\end{align*}%
The last line is precisely the result we claimed in lemma    \ref{lmm=controleOrliczEntropie}.

\section{A moment bound} %
\label{sec=controleDeMoments}
In this annex we prove lemma \ref{lmm=controleV}. The proof mainly follows the one in Miclo's doctoral dissertation, 
with a few changes to accomodate our hypotheses.
\subsection{Outline of the proof} %
\newcommand{\gen}{\mathcal{L}}
\newcommand{\genEps}{\gen_\epsilon}
\newcommand{\rhoEps}{\rho_\epsilon}
We need to introduce some notation.

For $\epsilon>0$, we denote by $\genEps$ the generator of the diffusion at fixed temperature $\epsilon$:
\[\genEps = \frac{\epsilon}{2}\Delta - \frac{1}{2} \grad V \grad\cdot.\]
We will need a smooth version of a step function; we call it $f$ and suppose that it satisfies:
\[%
  f(x) = \begin{cases}
    0 & \text{if } x\leq 0,\\
    \exp\left( - \exp\left(\frac{1}{x}\right)\right) & \text{on } [0,1],\\
    1 & \text {on } [2,\infty[.
  \end{cases}
\]
We recall the hypotheses on $V$:
\begin{itemize}
  \item It goes to infinity at infinity, 
  \item its gradient $\grad V$ is bounded, and 
  \item its Laplacian $\Delta V$ is negative for large $x$.
\end{itemize}
Note that, since $V$ is continuous, there must be an $R$ such that $\Delta V$ is negative whenever $V(x)\geq R$.

Finally, let $g$ be an increasing function, going to zero at zero.

The idea of the proof is that, as time goes by, the value of $V$ at $X_t$ has a typical scale, namely $\frac{1}{g(\sigma(t))}$, for a function $g$ to be made precise later, so that when we try to estimate $\esp(V^p(X_t))$, we only have to take into account the small values of $V$. 

More precisely, let $\rho_\epsilon(\cdot) =  f\left( g(\epsilon)V(\cdot) -
(R+1)\right)$. This is a smooth approximation of $\ind{V\geq
\frac{R}{g(\epsilon)}}$. We may bound the expectation of $V^p(X_t)$:
\begin{align}
 \notag \esp[V^p(X_t)] &= \esp[V^p\rho_{\sigma(t)}(X_t)] + \esp[V^p(1-\rho_{\sigma(t)}(X_t))] \\
\label{eq=majVp}
  &\leq \esp[V^p\rho_{\sigma(t)}(X_t)] + \left(\frac{R+3}{g(\sigma(t))}\right)^p.
\end{align}
To bound the first term, we use the explicit expression of the generator.
Intuitively, we write, for $h_t = V^p \rho_{\sigma(t)}$:
\[ \frac{d}{dt}(P_t h_t) = P_t \gen_{\sigma(t)}h_t + P_t (\frac{d}{dt} h_t), \]
and integrate between two times $t$ and $t'$. To ensure that everything exists, we use the stopping time
$T_k = \inf\{t, V(X_t) \geq k\}$. We get:
\newcommand{\tinfK}{{t\wedge T_k}}
\newcommand{\tpinfK}{{t'\wedge T_k}}
\[%
\begin{split} \esp[h_\tinfK (X_\tinfK)]%
  &= \esp [h_\tpinfK(X_\tpinfK)]\\
  &\quad + \esp[\int_{\tpinfK}^{\tinfK} \gen_{\sigma(s)}(h_s)(X_s)ds \\
  &\quad + \esp[\int_{\tpinfK}^{\tinfK} \sigma'(s)g'(\sigma(s))%
    f'\left(g(\sigma(s))V(X_s)-(R+1)\right)   V^{p+1}(X_s)  ds.
\end{split}
\]
Since $V$ is positive, $f$ and $g$ increasing and $\sigma$ decreases, the whole last term is negative. We try to estimate the second one, and study $\gen_{\sigma(s)}h_s(X_s)$.
\begin{lmm}\label{lmm=controleGen}
  Let us define $\varphi:x\mapsto x\log^2(x)$. There exists an $M$ and a time $t'$ (which may depend on $p$ and on the initial law) such that:
  \[  \forall t\geq t',\forall x, \quad %
  \gen_{\sigma(t)}(h_t)(X_t) \leq \exp\left(%
    - \frac{M}{\varphi\left(\sigma(t)g(\sigma(t))\right)}
    \right).
  \]
\end{lmm}
We postpone the proof and finish the argument. The inequality dictates the choice of $g$: $g = \ln(1/\cdot)^{-3}$ guarantees
\begin{align*}
\sigma(t)g(\sigma(t))& = \frac{1}{\ln(t)(\ln\ln(t))^3}, \\
\varphi(\sigma(t)g(\sigma(t))) &= \frac{\ln^{2}\left( 1/ \ln(t)(\ln\ln(t))^3\right) }{ \ln(t)(\ln\ln(t))^3 }
= \frac{\ln^2 \left( \ln(t)(\ln\ln(t))^3 \right)} {\ln(t) (\ln\ln(t))^3 }.
\end{align*}
Indeed, the upper bound on the generator then becomes
\begin{align*}
  \gen_{\sigma(t)}(h_t)(X_t) &\leq \exp\left(%
  - \frac{M}{\varphi(\sigma(t)g(\sigma(t))}
    \right) \\
    &\leq \exp\left(
      -M \ln(t)\times  \frac{(\ln\ln(t))^3}
            { \ln^2 \left( \ln(t)(\ln\ln(t))^3 \right) }
    \right).
\end{align*}
Since the ratio $ (\ln\ln(t))^3 / (\ln^2 (\ln(t) \ln\ln(t)^3))$ goes to infinity, it eventually exceeds $2/M$, so that for $t$ big enough, 
\[ 
  \gen_{\sigma(t)}(h_t)(X_t) \leq \exp\left(%
  - 2\ln(t) \right).
\]
  
Going back to the bound on the expected value we were looking for, the two previous arguments imply:
\[%
\esp[h_\tinfK(X_\tinfK)] \leq \esp[h_\tpinfK(X_\tpinfK)] 
  + \int_{t'}^{\infty}\exp\left(%
    - 2\ln(t) 
  \right).
\]
We succeeded in making the last integral finite. We can then let $K$ go to infinity, and since $t'$ is fixed, 
we get the existence of a constant $M_p$ (which depends on $p$ and on the initial law) such that:
\[%
\esp[h_t(X_t)] \leq M.
\]
Plugging this back into inequality \eqref{eq=majVp} yields:
\[
\esp[V^p(X_t)] \leq M + \left(\frac{R+3}{g(\sigma(t))}\right)^p.
\]
The expression of $g$ shows that, for a new constant $M$:
\[
  \esp[V^p(X_t)] \leq M (\sigma(t)\ln(t) (\ln\ln(t))^3)^p,
\]
and the result is proved.

\subsection{An estimate on the generator} %
We now turn to the proof of lemma \ref{lmm=controleGen}. We have to bound $\genEps(\rhoEps V^p)(x)$, and our first step 
will be to give a more explicit expression of this quantity. We will need the derivatives of $\rhoEps(x)$. To alleviate notations, 
we will write $y=y(x,\epsilon) = g(\epsilon)V(x) - (R+1)$.
\begin{align*}
  \rhoEps(x) &= f(g(\epsilon) V(x) - (R+1)) = f(y),\\
  \grad \rhoEps(x) &= g(\epsilon) f'(y)\grad V(x), \\
  \Delta \rhoEps(x) &= g(\epsilon)^2 f''(y) \norm{\grad V}^2 + g(\epsilon) f'(y) \Delta V.
\end{align*}
The quantity we would like to estimate is
\begin{align*}
  \genEps(\rhoEps V^p)(x)
  &= \rhoEps\genEps V^p(x) +  \epsilon\scal{\grad \rhoEps}{\grad V^p} (x) + V^p \genEps \rhoEps(x)
\end{align*}
We consider three cases, according to the value of $V(x)g(\epsilon)$.
\paragraph{$V$ is small: $V(x)g(\epsilon)\in[0,R+1]$} %
On this interval, $\rhoEps$ vanishes, so $\genEps(\rhoEps)$ is zero.

\paragraph{$V$ is large.} %
Let $\lambda$ be a strictly positive real, to be fixed later on. We consider the case where $V(x)g(\epsilon)\in [R+1+\lambda, \infty)$, which may be rewritten as: $y\in[\lambda,\infty)$. We develop the expression of $\genEps(\rhoEps V^p)$.
\[
\genEps(\rhoEps V^p)(x)%
  = \rhoEps\genEps V^p(x) +  \epsilon g(\epsilon) f'(y)\times p V^{p-1} \norm{\grad V}^2%
     + V^p \left( \frac{1}{2}\epsilon \Delta \rhoEps - \frac{1}{2}\scal{\grad \rhoEps}{\grad V}\right).
\]
We compute the derivatives of $\rhoEps$ and put together the terms involving $\norm{\grad V}^2$.
\begin{align*}
  \genEps(\rhoEps V^p)(x)%
  &= \rhoEps\genEps V^p(x) + \left( \epsilon g(\epsilon) f'(y)p V^{p-1} + V^p\left(%
    \frac{1}{2}\epsilon g(\epsilon)^2 f''(y) - \frac{1}{2}g(\epsilon) f'(y)\right)\right) \abs{\grad V}^2\\
    &\quad + \frac{1}{2}\epsilon g(\epsilon) f'(y) \Delta V\\
  &= A + B +C.
\end{align*}
Since $V\times g(\epsilon)\geq R$, $V\geq R$. We already noted that $R$ may be chosen so that, if $V$ is bigger than $R$, 
$\Delta V$ is less than zero, and this makes the third term $C$ negative. The term $B$ can be rewritten as:
\begin{align}
\notag%
B &= \left(\epsilon g(\epsilon) f'(y) pV^{p-1} + V^p\left(%
    \frac{1}{2}\epsilon g(\epsilon)^2 f''(y) - \frac{1}{2}g(\epsilon) f'(y)\right)\right) \abs{\grad V}^2\\
\label{equ=expressionDeB}%
  &= V^p g(\epsilon)\left( \left(\frac{p\epsilon}{V}- \frac{1}{2}\right)f'(y)%
     + \frac{1}{2}\epsilon g(\epsilon) f''(y)%
    \right)\abs{\grad V}^2.
\end{align}
We add another condition on $f$: it should be concave when $y$ is near $2$ (\eg{} on $[\frac{3}{2},2]$). 
On $[\lambda, 3\lambda/2]$, $f''/f'$ is bounded --- let $M$ be a bound. This entails:
\[\forall y\geq \lambda, \qquad f''(y) \leq M f'y).\]
Coming back to $B$, we deduce:
\[
B \leq \left(\frac{p\epsilon}{V} + \frac{M\epsilon g(\epsilon)}{2}%
- \frac{1}{2}\right) f'(y) g(\epsilon) V^p \abs{\grad V}^2.
\]
The term between brackets is negative, uniformly in $x$ as soon as $\epsilon$ is small enough.

Finally, the first term $A   = \genEps V^p$ is also negative:
\begin{align*}
A  &= \frac{\epsilon}{2}\Delta(V^p)  - \frac{1}{2}\scal{\grad V}{\grad(V^p)}\\
   &= \frac{\epsilon}{2}\left(%
     p(p-1)V^{p-1} \abs{\grad V}^2 + pV^{p-1} \Delta V
   \right)
     - \frac{1}{2}V^{p-1}\abs{\grad V}^2\\
   &\leq \left( \frac{p(p-1)\epsilon}{2} - \frac{1}{2}\right)\abs{\grad V}^2. 
\end{align*}
Once more, the term between brackets is negative when $\epsilon$ is small. 
To conclude, for any $\lambda$, there exists an $\epsilon_0$ such that:
\[%
  \forall \epsilon < \epsilon_0, \forall x,%
    \qquad V(x)g(\epsilon) \geq R+1+\lambda  \implies \genEps(\rhoEps V^p) \leq 0. 
\]

\paragraph{$V$ is of the order of $R/g(\epsilon)$.} %
This last case is that where $g(\epsilon) V(x) \in [R+1, R+1+\lambda]$. Let us reuse the decomposition $\genEps(\rhoEps V^p) = A+B+C$ from the previous paragraph. The same reasoning applies for $A$ and $C$, and they are both negative, so it suffices to get a bound on $B$. From \eqref{equ=expressionDeB}:
\[%
B = \left(%
      \left( \frac{p\epsilon}{V} - \frac{1}{2}\right) f'(y) + \frac{1}{2}\epsilon g(\epsilon) f''(y)
  \right)%
  g(\epsilon)V^p\abs{\grad V}^2.
\]
If we choose $R$ sufficiently big and $\epsilon$ small enough, the quantity between brackets in front of
$f'(y)$ is less than $1/4$.
\[%
  B \leq \left( -\frac{1}{4} f'(y) +\frac{1}{2} \epsilon g(\epsilon) f''(y) \right)%
  g(\epsilon)V^p \abs{\grad V}^2.
\]
Recall that    $f =  \exp(-\tau )$, where $\tau(y)= \exp(1/y)$. This implies:
\begin{align*}%
B &\leq \left( \frac{1}{4} \tau' f + \frac{1}{2}\epsilon g(\epsilon)(- \tau''f + (\tau')^2 f) \right) 
  g(\epsilon)V^p \abs{\grad V}^2\\
  &\leq \frac{1}{2}\left( \frac{1}{2}\tau' f + \epsilon g(\epsilon)(\tau'(y))^2 f(y)\right)
  g(\epsilon)V^p \abs{\grad V}^2\\
\end{align*}
Define $h_\epsilon = \frac{1}{2}\tau' f + \epsilon g(\epsilon)\tau'^2 f$. We study it by differentiating:
\[%
  h_\epsilon' = \left( \frac{1}{2}\tau'' - \frac{1}{2}\tau'^2 
    + 2\epsilon g(\epsilon) \tau'\tau'' - \epsilon g(\epsilon)\tau'^3\right) f.%
\]
The explicit expression of $\tau$ ensures:
\[\exists \lambda \forall y\in[0,\lambda] \qquad 0\leq \tau''(y)\leq \frac{1}{4}\tau'^2(y).\]
This $\lambda$ does not depend on $\epsilon$. This can be used to bound $h_\epsilon'$ from below:
\begin{align*}
  h_\epsilon'(y)
  &\geq \left( - \frac{1}{2}\tau'^2(y) + \frac{1}{2}\epsilon g(\epsilon) \tau'(y)^3%
     - \epsilon g(\epsilon) \tau'(y)^3\right) f(y)\\
  &\geq \left( - \frac{1}{2} - \frac{1}{2} \epsilon g(\epsilon) \tau'(y)\right)\tau'(y)^2 f(y).
\end{align*}
\newcommand{\yUnEps}{y_{1,\epsilon}}
\newcommand{\yDeuxEps}{y_{2,\epsilon}}
\newcommand{\yEps}{y_\epsilon}
Let $\yUnEps$ be the solution of the equation: $- 1 - \epsilon g(\epsilon) \tau'(y) = 0$. When $\epsilon$ is small, $\yUnEps$ 
will be less than $\lambda$, and the monotonicity of $\tau'$ will give:
\[%
\forall y \leq \yUnEps, \quad h_\epsilon'(y) \geq 0.
\]
Similarly, $h_\epsilon'$ can be bounded above:
\begin{align*}
  h_\epsilon' 
    &\leq \left( \frac{1}{8}\tau'^2(y) - \frac{1}{2}\tau'^2(y) - \epsilon g(\epsilon)\tau'^3(y)\right)f(y)\\
    &\leq \left( -\frac{3}{8} - \epsilon g(\epsilon) \tau'(y)\right)\tau'(y)^2 f(y).
\end{align*}
Now, let $\yDeuxEps$ be the root of $ - \frac{3}{8} - \epsilon g(\epsilon) \tau'(y) = 0$. Once more, when $\epsilon$ is small, $\yDeuxEps$ falls within $[0, \lambda]$. We deduce:
\[%
  \forall y\in [\yDeuxEps, \lambda], h_\epsilon'(y) \leq 0.
\]
We now know the $h_\epsilon$ increases on $[0,\yUnEps]$, and decreases on $[\yDeuxEps,\lambda]$, so that its maximum must be reached somewhere between these two points. More precisely, whenever $\epsilon$ is less than some $\epsilon_0$, it holds that
\[
  \exists \yEps\in [\yUnEps,\yDeuxEps],\forall y\in[0, \lambda], \quad h_\epsilon(y) \leq h_\epsilon(\yEps).
\]
The bounds on $\yEps$, the fact that $\tau$ decreases and the equations defining $\yUnEps, \yDeuxEps$ allow us to conclude:
\begin{align*}
  \forall y\leq \lambda, \qquad%
  h_\epsilon(y) &\leq \left(\frac{1}{2}\tau'(\yEps) + \epsilon g(\epsilon) \tau'(\yEps)^2\right) f(\yEps)\\
  &\leq \left( \frac{1}{2} \tau'(\yDeuxEps) + \epsilon g(\epsilon) \tau'(\yUnEps)^2\right)  f(\yDeuxEps)\\
  &\leq \left(  - \frac{3}{16\epsilon g(\epsilon)} + \frac{1}{\epsilon g(\epsilon)} \right) f(\yDeuxEps)\\
  &\leq \frac{1}{\epsilon g(\epsilon)}f(\yDeuxEps).
\end{align*}
It remains to estimate $f(\yDeuxEps) = \exp ( -\tau(\yDeuxEps))$. Since $\yDeuxEps$ is defined as a solution of an equation involving $\tau'$, we would like to compare $\tau$ and $\tau'$. The explicit expression of $\tau$ easily implies:
\[
 \ln(\abs{\tau'(y)}) = \ln(y^{-2}) + \frac{1}{y} \geq \frac{1}{y},
 \]
 therefore:
 \[
 \tau(y)  = y^2 \abs{\tau'(y)} \geq \frac{ \abs{\tau'(y)}}{\ln^2(\abs{\tau'(y)})}
 \]
 Applying this for $y=\yDeuxEps$, for which $\abs{\tau'(y)} = 3/(8 \epsilon g(\epsilon))$, entails:
 \begin{align*}
 \tau(\yDeuxEps)
 &\geq \frac{3}{8 \epsilon g(\epsilon) \ln^2 ( 8\epsilon g(\epsilon)/3)} \\
 & \geq\frac{3}{8 \epsilon g(\epsilon) \ln^2 (\epsilon g(\epsilon))}. 
 \end{align*}
Turning back to $f$, and defining $\varphi:x\mapsto x\ln^2(x)$, and $M = 3/8$, we have:
\begin{align*}
  f(\yDeuxEps)
  &= \exp( - \tau(\yDeuxEps)) 
  \leq \exp \left( - \frac{M}{\varphi(\epsilon g(\epsilon))} \right).\\
 \end{align*}
We now come back to the upper bound on $B$, and plug in the last equation.
\begin{align*}
  B &\leq \frac{1}{2}\times \frac{1}{\epsilon g(\epsilon)} \exp\left(%
  -\frac{M}{\varphi(\epsilon g(\epsilon))}
    \right)%
    g(\epsilon) V^p \abs{\grad V}^2.
\end{align*}
Since we suppose that $V(x)g(\epsilon)$ belongs to $[R+1,R+2]$, we may bound $V^p$ by $g(\epsilon)^{-p}$. We also supposed that $\grad V$ is bounded, so that there exists an $M'$ such that:
\[ B \leq \frac{M'}{\epsilon g(\epsilon)^p} \exp\left( - \frac{M}{\varphi(\epsilon g(\epsilon))}\right).\] 
Up to a slight change of the constant $M$ in the exponential, we may neglect the pre-exponential term and write:
\[
  B\leq M'' \exp\left(%
    - \frac{M}{\varphi(\epsilon g(\epsilon))}
    \right)
\]
This concludes the proof.

\section{Regularity results and estimates on the process}%
\subsection{An equivalent of the partition function} %
\label{sec=fonctionDePartition}
We recall here Laplace's method, which enable us to study the asymptotic behaviour of the partition function, \ie{} the constant $Z_\sigma = \int \exp ( -V/\sigma)dx$.

\begin{thm}%
  Let $V$ be a function from $\R^d$ to $\R$, satisfying hypotheses \ref{hpth=globalMinimum} and \ref{hpth=croissanceDeV} ($V$ has a unique, well behaved, global minimum, and $V$ goes to infinity at infinity rapidly enough).
  Then $Z_\sigma$ exists, and the following holds:
  \[ Z_\sigma \underset{\sigma\to 0}{\sim} \frac{1}{\sqrt{\det \hess V}} \left(\frac{\sigma}{2\pi}\right)^{d/2}. \]
\end{thm}%

To prove this classical result, we cut the integral in two parts, the main one (near the origin) and a remainder. Before we proceed, let us remark that, up to a change of coordinates, we may as well suppose that $\hess(V)_0$ is a diagonal matrix, and we have Taylor's formula:
\[ V(x) = \frac{1}{2}\sum_i \lambda_i x_i^2 + \epsilon(x)\sum_i x_i^2,\]
where $\epsilon(x)$ goes to zero at $0$. We choose an $r$ such that, on $B = [-r,r]^d$, $\epsilon(x)\leq \frac{1}{4}(\inf \lambda_i) \sum x_i^2$.

Let us begin by the negligible part, outside of $B$. Since $V$ goes to infinity, and $0$ is the unique global minimum, 
there exists an $\eta>0$ such that $V(x) \geq \eta$ outside $B$.
We introduce an $\exp(-V)$ in the integral (the growth hypothesis makes it integrable), and use this bound:
\begin{align*}%
  \int_{x \notin B} \exp ( -V/\sigma) dx
&= \int_{x\notin B} \exp(-V) \exp\left( -( 1/\sigma - 1)V(x) \right) dx \\
&\leq \int_{x\notin B} \exp(-V)dx \exp\left( -(1/\sigma - 1)\eta \right) \\
&\leq Z_1 \exp\left( -( 1/\sigma - 1)\eta\right).
\end{align*}%

Let us turn to the main term.  We divide it by $\sigma^{d/2}$ (so that we only have to find a limit). We change variables and use $x = \phi_\sigma(y)$ defined by $x_i = y_i \sqrt{\sigma/\lambda_i}$.
\begin{align*}%
  \sigma^{-d/2}
  \int_B \exp( -V/\sigma) dx_1\cdots dx_n
  &= \sigma^{-d/2} \int \ind{x\in B} \exp\left(
    -\frac{1}{2} \sum_i \frac{\lambda_i}{\sigma} x_i^2 
    + \epsilon(x) \sum_i x_i^2 \right)dx \\
  &= \frac{1}{\sqrt{\lambda_1\cdots\lambda_n}} 
  \int \ind{\phi_\sigma(y)\in B} \exp\left( - \frac{1}{2}\sum y_i^2 + \epsilon(\phi_\sigma(y)) \right)dy.
\end{align*}%
The function inside the integral converges pointwise to $\exp(-\sum y_i^2)$ when $\sigma$ goes to zero (because $\phi_\sigma(y)$ goes to zero for a fixed $y$). It is bounded from above by the integrable function $\exp( -\frac{1}{4} \sum y_i^2$, and we may apply Lebesgue's dominated convergence:
\[
  \sigma^{-d/2}
  \int_B \exp( -V/\sigma) dx_1\cdots dx_n
  \xrightarrow[\sigma\to 0]{} 
  \frac{(2\pi)^{d/2}}{\sqrt{\lambda_1\cdots \lambda_n}} 
\]
With the bound on the remainder, this gives the equivalent of $Z_\sigma$.

\subsection{Finiteness of the entropy and regularity} %
\label{sec=finiteEntropy}
\newcommand{\given}{|}
We begin by proving that the relative entropy $I_t$ is finite. %
To do this, we study directly the explicit density, which we know thanks to a Girsanov transform. We follow 
a proof from \cite{Roy99}, with a few minor changes to deal with the non-homogeneity in time.

Recall that the process $X$ is defined by the following SDE: 
\[ dX_t = \sqrt{\sigma(t)}dB_t - \frac{1}{2} \grad V(X_t) dt.\]
If we define a new reference martingale $M_t = \int_0^t \sqrt{\sigma(s)} dB_s$, we may define $X$ as the solution to the SDE: 
\[ dX_t = dM_t - \frac{1}{2} \grad V(X_t) dt.\]
Note that $M_t$ is just a Brownian motion under a (deterministic) change of time --- if we define $\tau(t) = \int_0^t \sigma(s) ds$, $M_{\tau^{-1}}(t)$ is a Brownian motion.
To find the density of the law of $X_t$ with respect to its equilibrium measure $\mu_t$, we decompose it in three terms:
\[%
  \frac{d\law(X_t)}{d\mu_t} = \frac{ d\law(X_t)}{d\law(M_t)} 
  \times \frac{d\law{M_t}}{d\lambda}
  \times \frac{d\lambda}{d\mu_t}.
\]%
To compute the first term, we use the (trajectorial) density of $X$ with respect to $M$, which is given by Girsanov's theorem:
\begin{align*}
  F &= \exp\left( -\frac{1}{2} \int \grad V(M_s) dM_s 
         - \frac{1}{2} \int_0^t \frac{\abs{\grad V}^2}{4}(M_s) d\bra{M}_s
       \right) \\
       &= \exp\left(  -\frac{1}{2} \int \grad V(M_s) dM_s 
           - \frac{1}{8} \int_0^t \abs{\grad V}^2(M_s) \sigma(s) ds
        \right).
\end{align*}
To get rid of the martingale term in the exponential, we apply Itô's formula to $V$ and the martingale $M$: 
\[ V(M_t) = V(x) + \int_0^t \grad V(M_s) dM_s + \frac{1}{2} \int_0^t \Delta V(M_s) d\bra{M}_s.\]
The functional $F$ may thus  be rewritten:
\[ F = \exp\left( 
        \frac{1}{2}V(x) - \frac{1}{2}V(M_t) + \int_0^t \left(
	    \frac{1}{4}\Delta V(M_s) - \frac{1}{8}\Delta V(M_s)
	  \right) \sigma(s) ds
	\right).
\]
The three densities we are looking for are:
\begin{align*}
  \frac{ d\law(X_t)}{d\law(M_t)}(M_t) &= f(M_t) =  \esp[ F \given \mathcal{F}_{\{t\}} ] \\
  \frac{d\law{M_t}}{d\lambda}(y) &= \exp ( - 2 v_t(y)) = (2\pi \tau(t))^{-d/2} \exp\left( - \frac{(x-y)^2}{2\tau(t)} \right) \\ 
  \frac{d\lambda}{d\mu_t}(y) &= Z_{\sigma(t)} \exp\left( - \frac{V(y)}{\sigma(t)} \right).
\end{align*}
We take the product of these terms; the last two quantities may be put into the conditional expectation, so that the density we are looking for (say $G$) may be written as:
\[
  G(M_t) = Z_{\sigma(t)} \esp \left[
      F \exp\left( \frac{V(M_t)}{\sigma(t)} - 2v_t(M_t) \right) \given \mathcal{F}_{\{t\}}
    \right]
\]
Let us now define $\gamma: x\mapsto x\log(x)$, and start to study $I_t$. By definition, $I_t = \int \gamma(G(y))d\mu_t(y)$. Since $G$ is best expressed as a conditional expectation, we rewrite $I_t$:
\begin{align}%
  I_t = \esp \left[ \gamma( G(M_t))  \frac{d\mu_t}{d\law(M_t)} \right] \notag \\
  \label{eq=ItCommeEsperance}
  &= \esp \left[ \gamma (G(M_t)) \frac{1}{Z_{\sigma(t)}} \exp \left(
        - \frac{V(M_t)}{\sigma(t)} + 2v_t(M_t)
      \right) \right].
\end{align}%
Since $\gamma$ is convex, we may apply Jensen's conditional inequality to $\gamma(G(M_t))$, and develop $\gamma$:
\begin{align*}
  \gamma(G(M_t)) &\leq \esp \left[  \gamma \left( 
      Z_{\sigma(t)} F \exp\left( \frac{V}{\sigma} - 2v_t 
    \right)\right) \given \mathcal{F}_{\{t\}}
   \right] \\
   &\leq \esp\left[
     Z_\sigma F \exp\left( \frac{V}{\sigma} - 2v_t \right) \left( \log Z_\sigma + \log F + \frac{V}{\sigma} - 2 v_t\right)
     \given \mathcal{F}_{\{t\}}
     \right].
\end{align*}
Multiply both sides by $\exp( -V/\sigma + 2v_t)$, and take the expected value; the left hand side becomes $I_t$ (thanks to \eqref{eq=ItCommeEsperance}), the conditioning disappears and we get:
\[
I_t \leq \esp \left[ F \left(
                 \log Z_\sigma + \log F + \frac{V(M_t)}{\sigma(t)} - 2v_t(M_t)
	    \right) \right]
\]
Recall that $F$ is a density, so that $\esp[F] = 1$, and we may take the constant $Z_\sigma$ out of the expectation. We add and substract $(2/\sigma)\log(F)$ inside the integral --- this will help us get rid of the term $V(M_t)/\sigma$:
\[
I_t \leq \log(Z_\sigma) - (2/\sigma - 1) \esp\left[ F \log F \right]
 + \esp\left[ F\left( \frac{2}{\sigma} \log F + \frac{V(M_t)}{\sigma} - 2v_t(M_t)\right) \right].
 \]
 Since $x\log x$ is bounded below, and $2/\sigma -1$ is positive, the second term is bounded from above (for any finite time $t$). The same is true for the first term. The only thing to check is that the last term is finite; let us call this term $A$. Since $F$ is given by an exponential, $A$ is given by:
 \[ A = \esp\left[ F\left( \frac{1}{\sigma} V(x) + \frac{1}{4\sigma(t)}\int_0^t \left(
     2 \Delta V - \abs{\grad V}^2 \right)(M_s) \sigma(s) ds 
      - 2v_t(M_t) \right) \right].
\]
Let us consider the quantity between brackets. The first term is finite and does not depend on $M_t$. The integral is bounded above by something also independant of $M_t$ (indeed, $2\Delta V - \abs{\grad V}^2$ is uniformly bounded from above, because $\Delta V$ is negative outside a compact set).
The only thing left to check is that:
\[ \esp \left[ F ( - 2 v_t(M_t)) \right] <\infty. \]
We have already seen the explicit value of $v_t$: 
\[ \exp\left( -2 v_t(y) \right) = ( 2\pi \tau(t))^{-d/2} \exp\left( - \frac{(y-x)^2}{2\tau(t)} \right).\]
Taking logarithms, we see that:
\[ -2v_t(y) = -\frac{d}{2}\log(2\pi\tau(t)) - \frac{(y-x)^2}{2\tau(t)}.\]
Since the last term is positive, this quantity is bounded from above by something which does not depend on $y$. Therefore, $\esp[ -F\times(2v_t(M_t))]$ is finite. This concludes the proof.

\tableofcontents
\bibliographystyle{smfalpha}
\bibliography{biblio}

\end{document}